\documentclass[preprint,reqno]{imsart}

\usepackage{amsthm,amsmath,natbib}
\RequirePackage[colorlinks,citecolor=blue,urlcolor=blue]{hyperref}
\usepackage{amsfonts}
\usepackage{graphicx}
\usepackage{comment}
\usepackage[usenames,dvipsnames]{xcolor}
\usepackage{enumerate}
\usepackage{extarrows}
\usepackage{hyperref}
\usepackage{stackrel}
\usepackage{centernot}
\usepackage{caption}
\usepackage{subcaption}
\usepackage{fancyhdr}
\usepackage{tikz}
\usepackage{float}
\arxiv{????}

\startlocaldefs

\def\boxit#1{\vbox{\hrule\hbox{\vrule\kern6pt
          \vbox{\kern6pt#1\kern6pt}\kern6pt\vrule}\hrule}}

\newtheorem{theorem}{Theorem}
\newtheorem{corollary}{Corollary}

\newtheorem{lemma}{Lemma}

\numberwithin{equation}{section}
\numberwithin{theorem}{section}
\numberwithin{corollary}{section}
\numberwithin{proposition}{section}
\numberwithin{lemma}{section}
\numberwithin{definition}{section}
\numberwithin{remark}{section}

\newcommand{\veps}{\varepsilon}

\newcommand{\todr}{\stackrel{\mathrm{D}}{\longrightarrow}}
\newcommand{\eqdr}{\stackrel{\mathrm{D}}{=}}

\newcommand{\R}{\Bbb{R}}

\newcommand{\N}{\Bbb{N}}

\newcommand{\rmd}{{\rm d}}
\newcommand{\rmi}{{\rm i}}
\newcommand{\halmos}{\quad\hfill\mbox{$\Box$}}

\newcommand{\wh}{\widehat}
\newcommand{\what}{\widehat}

\newcommand{\dto}{\downarrow}

\newcommand{\be}{\begin{equation}}
\newcommand{\ee}{\end{equation}}
\newcommand{\bea}{\begin{eqnarray}}
\newcommand{\eea}{\end{eqnarray}}
\newcommand{\bean}{\begin{eqnarray*}}
\newcommand{\eean}{\end{eqnarray*}}
\newcommand{\ben}{\begin{equation*}}
\newcommand{\een}{\end{equation*}}
\newcommand{\ba}{\begin{aligned}}
\newcommand{\ea}{\end{aligned}}

\def\nexto{\kern -0.54em}

\newcommand{\PP}{\textbf{\rm P}}
\newcommand{\EE}{\textbf{\rm E}}
\newcommand{\PD}{\textbf{\rm PD}}

\newcommand{\bfa}{{\bf a}}
\newcommand{\bfy}{{\bf y}}
\newcommand{\bfA}{{\bf A}}
\newcommand{\bfD}{{\bf D}}

\newcommand{\bfM}{{\bf M}}
\newcommand{\bfm}{{\bf m}}
\newcommand{\bfp}{{\bf p}}
\newcommand{\bfQ}{{\bf Q}}

\newcommand{\bfMult}{{\bf Mult}}

\newcommand{\lf}{\lfloor}
\newcommand{\rf}{\rfloor}
\newcommand{\lc}{\lceil}
\newcommand{\rc}{\rceil}

\endlocaldefs

\begin{document}

\begin{frontmatter}

\title{Limiting Behaviour of  Poisson-Dirichlet and Generalised  Poisson-Dirichlet Distributions}
\runtitle{Limiting Behaviour of  Poisson-Dirichlet distributions }

\author{Ross Maller$^a$ and Soudabeh Shemehsavar$^{c}$
\thanks{$^a$Research School of Finance, Actuarial Studies \& Statistics, 
Australian National University, Ross.Maller@anu.edu.au;  
\newline
$^c$College of Science, Technology, Engineering and Mathematics, Murdoch University, and
School of Mathematics, Statistics \& Computer Sciences, University of Tehran, 
\newline
shemehsavar@ut.ac.ir (corresponding author)
}}

\begin{abstract}
We derive large-sample and other limiting distributions of  the 
``frequency of frequencies'' vector, $\bfM_n$, together with the number of species,   $K_n$, in a  Poisson-Dirichlet or generalised  Poisson-Dirichlet gene or species sampling model.
Models analysed include those  constructed from gamma and  $\alpha$-stable subordinators by Kingman, the  two parameter extension by Pitman and Yor, and another  two parameter version constructed by omitting large jumps from  an $\alpha$-stable subordinator.
In the Poisson-Dirichlet case $\bfM_n$ and $K_n$ turn out to be  asymptotically independent, and notable, especially for statistical applications, is that in other cases 
 the conditional limiting distribution of $\bfM_n$, given $K_n$, is normal, after certain centering and norming.
\end{abstract}



\end{frontmatter}

\section{Introduction}\label{intro}
\cite{kingman1975} suggested a way of constructing random distributions on the infinite unit  simplex
by ranking the  jumps  of a subordinator up to a specified 
time,
then normalising them by the value of the subordinator at that time.
Taking the subordinator to be a driftless  gamma subordinator 
generates his well known Poisson-Dirichlet distribution $\PD(\theta)$, 
for $\theta>0$,
which was later shown to be intimately connected with the famous Ewens sampling formula in genetics. 
Another of Kingman's distributions,  denoted by $\PD_\alpha$,  arises when a driftless $\alpha$-stable subordinator with parameter $\alpha\in(0,1)$ is used instead of the gamma subordinator. Later again, an encompassing two parameter  Poisson-Dirichlet distribution $\PD(\alpha, \theta)$ was constructed by 
\cite{PY1997}.

These distributions and the methodologies associated with them have subsequently had a big impact in many applications areas,  especially in population genetics, but also in the excursion theory of stochastic processes,  the theory of  random partitions, random graphs and networks, probabilistic number theory,
machine learning,  Bayesian statistics, and a variety of others.
They have also given rise to a number of generalisations and a large literature analysing the various versions. 
Among these, generalising Kingman's $\PD_\alpha$ class,
 is the two parameter  $\PD_\alpha^{(r)}$ class of
\cite{IpsenMaller2017}, defined for $\alpha\in(0,1)$ and $r>0$.

We consider a sample of size $n\in \N=\{1,2,\ldots\}$ from  $\PD_\alpha$,  $\PD(\theta)$,  $\PD(\alpha,\theta)$
 or  $\PD_\alpha^{(r)}$, 
and form the ``frequency of frequencies'' random vector, $\bfM_n$, together with the number of species,  random variable $K_n$.
Thus we have $\bfM_n= (M_{1n}, M_{2n},  \ldots, M_{nn})$, where $M_{jn}$ is the number of alleles or species represented $j$ times in the sample,  and 
$K_n =M_{1n}+M_{2n}+ \cdots+ M_{nn}$ is the total number of alleles or species in the sample. 
Specialising  our general notation, when the model under consideration depends on parameters,  these are made explicit in the particular 
$\bfM_n$ and $K_n$ analysed; thus, from $\PD_\alpha$ we have $(\bfM_n(\alpha),K_n(\alpha))$,
from $\PD(\theta)$ we have $(\bfM_n(\theta),K_n(\theta))$, etc.

The paper is organised as follows. We start in the next section with 
$\PD_\alpha^{(r)}$, finding the limiting distribution of 
$(\bfM_n(\alpha,r),\, K_n(\alpha,r))$
as $n\to\infty$. In  Section \ref{s3} we apply the methods thus developed to
 $\PD_\alpha$,  while Section \ref{s4} has $\PD(\theta)$   and $\PD(\alpha,\theta)$.
 Section \ref{s5} gives other related limiting results;
for example, the limit as $r\dto 0$ of
 $\PD_\alpha^{(r)}$ is $\PD_\alpha$.
     Special cases of some of the results are known of course but for each  we  add some extra information.
    In particular, we obtain the joint limiting distribution of the appropriate  $(\bfM_n,K_n)$, in each case.
    For $\PD(\theta)$, it   turns out that
 $\bfM_n(\theta)$ and $K_n(\theta)$ are asymptotically independent, while for  $\PD_\alpha^{(r)}$, $\PD_\alpha$ and $\PD(\alpha,\theta)$,
 the conditional  limiting distribution of $\bfM_n$ given  $K_n$ is normal, after certain centering and norming --
 a useful fact for statistical applications.
 Further discussion of these results and more references are in the concluding Section \ref{s6}.

\section{Limiting Distribution of $(\bfM_n(\alpha,r),\, K_n(\alpha,r))$  as $n\to\infty$}\label{s2}
In 
 the $\PD_\alpha^{(r)}$ model $\bfM_n(\alpha,r)$ and $K_n(\alpha,r)$ depend on the parameters $\alpha\in(0,1)$ and $r>0$, and the sample size 
 $n\in \N$.
In this section 
$\alpha$ and $r$ are kept fixed for the large-sample analysis ($n\to\infty$).
From Eq. (5.11) of    \cite{IpsenMallerShemehsavar2021} 
 (hereafter, ``IMS'')
we get the following  formula for the distribution of the ($n$+1)--vector
 $(\bfM_n(\alpha,r),K_n(\alpha,r))$:
\be\label{4.1}
\PP(\bfM_n(\alpha,r)=\bfm,\, K_n(\alpha,r)=k)
= n
\int_{0}^{\infty}\frac{\Gamma(r+k)\lambda^{\alpha k}}
{\Gamma(r)\Psi(\lambda)^{r+k}}
\prod_{j=1}^n \frac{1}{m_j!}\big(F_j(\lambda)\big)^{m_j}\frac{\rmd \lambda}{ \lambda}
\ee
where
\be\label{2.7aa}
\Psi(\lambda)=1+\alpha \int_{0}^1(1-e^{-\lambda z})z^{-\alpha-1} \rmd z
\ee
and
\be\label{Fdef}
F_j(\lambda) :=\frac{\alpha}{j!} \int_0^\lambda e^{-z}  z^{j-\alpha-1}\rmd z, 
\ j\in\N_n=\{1,2,\ldots,n\},\ \lambda>0.\\
\ee
In \eqref{4.1}, $k\in \N_n$, $n\in\N$, and $\bfm$ takes values in 
\be\label{4.4a0}
A_{kn}:=\Big\{\bfm= (m_1,\ldots,m_n): m_j\ge 0, \, \sum_{j=1}^njm_j=n,\, \sum_{j=1}^nm_j=k\Big\}.
\ee

Introduce for each $\lambda>0$ a subordinator 
$(Y_t(\lambda))_{t>0}$ having 
L\'evy measure 
\be\label{LYP}
 \Pi_\lambda(\rmd y) :=
\frac{ \alpha y^{-\alpha-1}\rmd y} {\Gamma(1-\alpha)}
 \Big(
 {\bf 1}_{\{0<y<\lambda\le 1\}} +  {\bf 1}_{\{0<y<1<\lambda\}}\Big).
 \ee
As shown in IMS,  each $Y_t(\lambda)$, $t>0$, $\lambda>0$,  has a continuous bounded density which we denote by
$f_{Y_t(\lambda)}(y)$, $y>0$.
Let $J\ge 1$ be a fixed integer,  define 
\be\label{qdef}
q_j = \frac{\alpha \Gamma(j-\alpha)}{j!\Gamma(1-\alpha)}, \ j\in \N_J, 
\ee
let $Q_1=1$, and  when $J\ge 2$, set
$Q_j=(1-\sum_{i=1, i\ne j}^Jq_i) \prod_{i=1, i\ne j}^J  q_i$, $2\le j\le J$.
Let
 $Q= \prod_{j=1}^Jq_j$, 
and let $\bfQ$ be the $J\times J$ matrix whose  inverse is
\be\label{Qdef} 
\bfQ^{-1}=
\frac{1}{(1-\sum_{j=1}^Jq_j)\prod_{j=1}^Jq_j}
 \begin{bmatrix}
   Q_1 & Q & \cdots & Q \\
   Q & Q_2 & \cdots & Q \\
   \vdots  & \vdots  & \ddots & \vdots  \\
  Q & Q & \cdots & Q_J 
 \end{bmatrix}.
\end{equation}
Let   $\bfa=(a_1,a_2, \ldots, a_J)\in\R^J$, $c>0$,
$m_j\ge 0$,  $1\le j\le n$, 
and  for the components of $\bfM_n(\alpha,r)$ write $(M_{jn}(\alpha,r))_{1\le j\le n}$.
Then we get:

\begin{theorem}\label{M4}
For the $\PD_\alpha^{(r)}$ model, we have the  asymptotic:
\bea\label{30}
&&
\lim_{n\to\infty}
 \PP\Big(
\frac{M_{jn}(\alpha,r)}{K_n(\alpha,r)}\le 
q_j+ \frac{a_j}{n^{\alpha/2}},\, j\in\N_J,
\frac{K_n(\alpha,r)}{n^\alpha} \le c\Big)  \cr
&&\cr
&& = 
 \frac{1}{\Gamma(r)\Gamma^r(1-\alpha)}\times
  \int_{\bfy\le \bfa}
\int_{x\le c} \int_{\lambda>0}
\frac{x^{r+J/2-1} e^{-\tfrac{x}{2} \bfy^T\bfQ^{-1}\bf y}}
{\sqrt{(2\pi)^{J}{\rm det}(\bfQ)}}\cr
&&\cr
&& \hskip2cm
\times
e^{-x(\lambda^{-\alpha}\vee 1)/\Gamma(1-\alpha)}
f_{Y_x(\lambda)}(1)
\frac{\rmd \lambda}{ \lambda^{\alpha r+1}} \rmd x\rmd \bfy,
\eea
where $\bfy =(y_1,y_2,\ldots,y_J)\in\R^J$. 
The limiting distribution of 
\ben
n^{\alpha/2}\Big(
\frac{M_{jn}(\alpha,r)}{K_n(\alpha,r)}-q_j\Big)_{j\in\N_J},
\een
conditional on $K_n(\alpha,r)=x>0$, is $N(\bf0,\bfQ/{\rm x})$,
that is, with density
\be\label{2.10}
\frac{x^{J/2}e^{-\tfrac{x}{2} \bfy^T\bfQ^{-1}\bf y}}
{\sqrt{(2\pi)^{J}{\rm det}(\bfQ)}}, \ \bfy\in \R^J.
\ee
\end{theorem}
\noindent {\bf Remark.}\
Integrating out $\bfy$ from \eqref{30} gives 
\begin{align}\label{M30}
&
 \frac{x^{r-1} }{\Gamma(r)\Gamma^r(1-\alpha)}
\int_{\lambda>0}
e^{-x(\lambda^{-\alpha}\vee 1)/\Gamma(1-\alpha)}
f_{Y_x(\lambda)}(1)
\frac{\rmd \lambda}{ \lambda^{\alpha r+1}}, \ x>0,
\end{align} 
for the limiting density of 
$K_n(\alpha,r)/n^\alpha$ at $x>0$, agreeing with (2.8) of  IMS.

\medskip
\noindent {\bf Proof of Theorem \ref{M4}:}\
Taking $n>J$, from \eqref{4.1} we can write
\bea\label{M4.1}
&&
\PP(M_{jn}(\alpha,r)=m_j, j\in\N_J,\, K_n(\alpha,r)=k)=\cr
&&\cr
&&
n
\int_{0}^{\infty}\frac{\Gamma(r+k)\lambda^{\alpha k}}
{\Gamma(r)\Psi(\lambda)^{r+k}}
\prod_{j=1}^J \frac{1}{m_j!} \big(F_j(\lambda)\big)^{m_j}\sum_{\bfm^{(J)}\in A_{kn}^{(J)}}
\prod_{j=J+1}^n \frac{1}{m_j!} \big(F_j(\lambda)\big)^{m_j}\frac{\rmd \lambda}{ \lambda}\cr
&&
\eea
where $\bfm^{(J)}=(m_{J+1}, \ldots, m_n)$
and
\be\label{A2def} 
A_{kn}^{(J)}=\Big\{m_j\ge 0, J< j\le n: \,
 \sum_{j=J+1}^njm_j=n-m_{++},\, 
\sum_{j=J+1}^n m_j=k- m_{+}\Big\},
\ee
with $m_{+}= \sum_{j=1}^Jm_j$ and
$m_{++}= \sum_{j=1}^Jjm_j$.
For each $\lambda>0$ let
\be\label{M3}
\bfp_n^{(J)}(\lambda)= \big(p_{jn}^{(J)}(\lambda)\big)_{J+1\le j\le n} 
=
\left(\frac{F_j(\lambda)}
{\sum_{\ell=J+1}^n F_\ell(\lambda)}\right)_{J+1\le j\le n}
\ee
and let 
 $\bfMult_{kn}^{(J)}(\lambda)$ 
 be a multinomial $(k-m_{+},\bfp_n^{(J)}(\lambda))$ vector with
 distribution 
\be\label{M0}
\PP\big(\bfMult_{kn}^{(J)}(\lambda)=(m_{J+1},\ldots, m_n) \big)=
(k-m_{+})!
\prod_{j=J+1}^n
\frac{\big(p_{jn}^{(J)}(\lambda)\big)^{m_j}}{m_j!},
\ee
 where $m_j\ge 0$,  $J+1\le j\le n$, and $\sum_{j=J+1}^n m_j=k-m_{+}$.
Just as in  IMS  the 
summation over $\bfm^{(J)}$ of the LHS 
of \eqref{M0} can be written as
\be\label{M9}
\PP\Big(\sum_{i=1}^{k-m_{+}} X_{in}^{(J)}(\lambda)=n- m_{++}\Big),
 \ee
 where $\big(X_{in}^{(J)}(\lambda)\big)_{1\le i\le k-m_{+}}$ are i.i.d. with
\be\label{M8}
\PP\big(X_{1n}^{(J)}(\lambda)=j\big) =p_{jn}^{(J)}(\lambda),\ J+1\le j\le n.
\ee
For brevity let $k'=k-m_{+}$ and $n'=n- m_{++}$ where convenient in what follows.
Then we have
\bea\label{M01}
&&\sum_{\bfm^{(J)}\in A_{kn}^{(J)}}
\prod_{j=J+1}^n \frac{1}{m_j!} \big(F_j(\lambda)\big)^{m_j} \cr
&&\cr
&&=
\frac{1}{k'!}
\sum_{\bfm^{(J)}\in A_{kn}^{(J)}}
k'!
\prod_{j=J+1}^n \frac{1}{m_j!} \big(p_{jn}^{(J)}(\lambda)\big)^{m_j}
 \Big(\sum_{\ell=J+1}^n F_\ell(\lambda)\Big)^{m_j}\cr
 &&\cr
 &&=
\frac{1}{ k'!} \Big(\sum_{\ell=J+1}^n F_\ell(\lambda)\Big)^{k'}
\PP\Big(\sum_{i=1}^{k'} X_{in}^{(J)}(\lambda)=n'\Big).
\eea
So we can write, from \eqref{M4.1},
\bea\label{M4.2aa}
&&
\PP\big(M_{jn}(\alpha,r)=m_j, \,
j\in\N_J,\,  K_n(\alpha,r)=k\big)\cr
&&\cr
&&= n
\int_{0}^{\infty}\frac{\Gamma(r+k)\lambda^{\alpha k}}
{\Gamma(r)\Psi(\lambda)^{r+k}}
\prod_{j=1}^J \frac{1}{m_j!} \big(F_j(\lambda)\big)^{m_j}\cr
&& \times
\frac{1}{ k'!} \Big(\sum_{\ell=J+1}^n F_\ell(\lambda)\Big)^{k'}
\PP\Big(\sum_{i=1}^{k'} X_{in}^{(J)}(\lambda)=n'\Big)
\frac{\rmd \lambda}{ \lambda}.
\eea
In this we change variable from $\lambda$ to $\lambda n$, which just means replacing $\lambda$ by $\lambda n$ in the functions in the integrand.
Then letting
\be\label{defqL}
q_{jn}=\frac{F_j(\lambda n)}
{\sum_{\ell=1}^n F_\ell(\lambda n)}, \ j\in\N_J,
\ {\rm and}\
q_{+n}= \sum_{j=1}^J q_{jn},
\ee
 so that
 \ben
\frac{\sum_{\ell=J+1}^n F_\ell(\lambda n)}
{\sum_{\ell=1}^n F_\ell(\lambda n)}=
1-\sum_{j=1}^J q_{jn}
=1- q_{+n}, 
\een
we can rewrite the RHS of \eqref{M4.2aa} as
\bea\label{M4.2}
&&
\int_{0}^{\infty}
\frac{\Gamma(r+k)}{k'!\prod_{j=1}^Jm_j!}
(1-q_{+n})^{k- m_{+}}\prod_{j=1}^Jq_{jn}^{m_j}
 \times
\frac{(\lambda n) ^{\alpha k}}{\Psi(\lambda n)^{k}}
\Big(\sum_{\ell=1}^n F_\ell(\lambda n)\Big)^{k}
\cr
&&\hskip2.5cm
\times
\frac{1}{\Gamma(r)\Psi(\lambda n)^{r}}
\times
n\PP\Big(\sum_{i=1}^{k'} X_{in}^{(J)}(\lambda n)=n'\Big)
\frac{\rmd \lambda}{ \lambda}.
\eea
Hence we have
\bea\label{M30a}
&&
 \PP\Big(
\frac{M_{jn}(\alpha,r)}{K_n(\alpha,r)}\le 
q_j+a_j/n^{\alpha/2},\, j\in\N_J,
\frac{K_n(\alpha,r)}{n^\alpha} \le c\Big)\cr
&&\cr
&&= 
\sum_{\bfm^{(J)},k}
{\bf 1}\big\{m_j \le (q_j+a_j/n^{\alpha/2})k,\, j\in\N_J,
k\le cn^\alpha\big\}\cr
&&\hskip1.5cm  \times
\PP(M_{jn}(\alpha,r)=m_j, \, j\in\N_J, K_n(\alpha,r)=k)\cr
&&\cr
&&= \int_{\bfm^{(J)},z} 
{\bf 1}\big\{m_j \le (q_j+a_j/n^{\alpha/2})z,\, j\in\N_J,
z\le cn^\alpha\big\}\cr
&&\cr
&&  \times
\PP(M_{jn}(\alpha,r)=\lf m_j\rf,\, j\in\N_J, K_n(\alpha,r)=\lf z \rf)
\prod_{j=1}^J \rmd m_j \rmd z,
\eea
where now the $m_j$ are treated as continuous variables.
Make the change of variables $z=n^\alpha x$,
$m_j= (q_j+ y_j/n^{\alpha/2})z$, 
so  $\rmd z=n^\alpha\rmd x$,
$\rmd m_j= \rmd y_j xn^{\alpha/2}$, 
to get the last expression as
\bea\label{M31}
&&
n^{\alpha(1+J/2)}
 \int_{\bfy,x} x^J
\PP\big(M_{jn}(\alpha,r)=
\lf (q_j+ y_j/n^{\alpha/2})xn^\alpha\rf, j\in\N_J, \cr
&&\cr
&& \hskip1.5cm  \times 
K_n(\alpha,r)=\lf xn^\alpha \rf\big)
\rmd \bfy\rmd x.
\eea
To find the limit as $n\to\infty$ of the probability in this expression, 
we use  \eqref{M4.2}.
The limits of the four factors in \eqref{M4.2}
are given in the next lemma, 
whose  proof  is deferred to an appendix.

\begin{lemma}\label{lemm1}
With the substitutions $k=\lf xn^\alpha \rf$, 
$m_j= (q_j+ y_j/n^{\alpha/2})k$, 
$k'=k- m_{+}$,
we have the following limiting behaviours as $n\to\infty$:
\be\label{L1}
\frac{\Gamma(r+k)}{k'!\prod_{j=1}^Jm_j!}
(1-q_{+n})^{k'}\prod_{j=1}^Jq_{jn}^{m_j}
\sim
\frac{(xn^\alpha)^{r-J/2-1}e^{-\tfrac{x}{2} \bfy^T\bfQ^{-1}\bf y}}
{ \sqrt{(2\pi)^{J}{\rm det}(\bfQ)}};
\ee
\bea\label{L2}
\lim_{n\to\infty}
\frac{(\lambda n)^{\alpha k}}{\Psi(\lambda n)^{k}}
\Big(\sum_{\ell=1}^n F_\ell(\lambda n)\Big)^{k}
=
e^{-x (\lambda^{-\alpha}\vee 1)/
\Gamma(1-\alpha)};
\eea
\bea\label{L3}
\frac{1}{\Gamma(r)\Psi(\lambda n)^{r}}
\sim
\frac{1}{\Gamma(r)(\lambda n)^{r\alpha}\Gamma^r(1-\alpha)};
\eea
\be\label{L4}
\lim_{n\to\infty} 
n\PP\Big(\sum_{i=1}^{k-m_{+}} X_{in}^{(J)}(\lambda n)=n- m_{++}\Big)= f_{Y_x(\lambda)}(1).
\ee
\end{lemma}

Multiplying together the right-hand sides of \eqref{L1}--\eqref{L4}, and keeping in mind the factor of $n^{\alpha(1+J/2)}x^J$ from \eqref{M31}, then substituting in \eqref{M4.2aa} and \eqref{M30a}, 
gives \eqref{30}, subject to verifying the various interchanges of integration and limit.
This can be done, but is unnecessary:
 \eqref{M30} specifies a proper probability distribution (total mass 1), 
 because we showed in  IMS that the limiting distribution of $K_n(\alpha,r)/n^\alpha$ is proper, and so the same is true after multiplying in the limiting normal distribution for the frequency spectrum.
This validates  the  convergence in \eqref{30}, and the 
 limiting conditional density  follows immediately from it.
 \halmos

 \section{Large sample limit for $\PD_\alpha$}\label{s3}
 In this section we obtain the large sample distribution for $\PD_\alpha$.
 It is known that for this, $K_n(\alpha)/n^\alpha$ has a limiting Mittag-Leffler distribution
 (Thm. 3.8, p.68,  of \cite{Pitman2006}).
We extend this to obtain the joint limit of   $\bfM_n(\alpha)$ with $K_n(\alpha)$.
An interesting part of the proof of the next theorem is  how  the density $f_{Y_x(\lambda)}(y)$ connects with the Mittag-Leffler density.
 Here a result of  \cite{Covo2009a}  plays a key role.
 Surprising also is the $n^{\alpha/2}$ norming for the frequency spectrum, as also occurred in Theorem \ref{M4}.
    Define $f_{Y_x(\lambda)}(y)$, $q_j$ and $\bfQ$ as in Section \ref{s2}
  (see \eqref{LYP},  \eqref{qdef} and \eqref{Qdef}),
 let $\bfa=(a_1,a_2, \ldots,a_J)\in\R^J$, $c\ge 0$, 
 $\bfy =(y_1,y_2,\dots,y_J)\in\R^J$
 and write $(M_{jn}(\alpha))_{1\le j\le n}$ for the components of 
 $\bfM_n(\alpha)$.
 
\begin{theorem}\label{M5}
For the $\PD_\alpha$ model we have the  asymptotic property:
\bea\label{330}
&&
\lim_{n\to\infty}
 \PP\Big(
\frac{M_{jn}(\alpha)}{K_n(\alpha)}
\le 
q_j+ \frac{a_j}{n^{\alpha/2}},\, j\in\N_J,
\frac{K_n(\alpha)}{n^\alpha} \le c\Big)
\cr
&&\cr
&& 
=
  \int_{\bfy\le \bfa} \int_{0<x\le c} 
\frac{x^{J/2-1}e^{-\tfrac{x}{2} \bfy^T\bfQ^{-1}\bf y}}
{\alpha  \sqrt{(2\pi)^J{\rm det}(\bfQ)}}
e^{-x/\Gamma(1-\alpha)}
f_{Y_x(1)}(1)     \rmd x\, \rmd \bfy.
\eea
Integrating out $\bfy$ gives 
\be\label{3M30}
 \frac{1}{\alpha x} 
e^{-x/\Gamma(1-\alpha)}
f_{Y_x}(1), \ x>0,
\ee
for the limiting density of 
$K_n(\alpha)/n^\alpha$.
This can alternatively be written as  a Mittag-Leffler density,
$f_{L_\alpha}(x)$, $x>0$.
The limiting density  of 
\ben
n^{\alpha/2}\Big(
\frac{M_{jn}(\alpha)}{K_n(\alpha)}-q_j\Big)_{j\in\N_J},
\een
conditional on $K_n(\alpha)=x>0$, is
$N(\bf0,\bfQ/{\rm x})$, as in \eqref{2.10}.
\end{theorem}
    
\noindent {\bf Proof of Theorem \ref{M5}:}\
This is a modification, in fact simplification, of the proof of Theorem \ref{M4}, so we  give just the main steps.
Starting from
\be\label{3.15}
\PP(\bfM_n(\alpha)=\bfm,\, K_n(\alpha)=k)
= \frac{n(k-1)!}{\alpha}
\Big(\frac{\alpha}{\Gamma(1-\alpha)}\Big)^k
\prod_{j=1}^n 
\frac{1}{m_j!} 
\Big(\frac{\Gamma(j-\alpha)}{j!}\Big)^{m_j}
\ee
(\cite{Pitman2006}, p.61), 
we get
\bea\label{M3.1}
&&
\PP(M_{jn}(\alpha)=m_j, j\in\N_J,\, K_n(\alpha)=k)\cr
&&\cr
&&
=
 \frac{n(k-1)!}{\alpha}
\Big(\frac{\alpha}{\Gamma(1-\alpha)}\Big)^k
\prod_{j=1}^J\frac{1}{m_j!} 
\Big(\frac{\Gamma(j-\alpha)}{j!}\Big)^{m_j}\cr
&&\cr
&& \hskip3.5cm  
\times
\sum_{\bfm^{(J)}\in A_{kn}^{(J)}}
\prod_{j=J+1}^n
\frac{1}{m_j!} 
\Big(\frac{\Gamma(j-\alpha)}{j!}\Big)^{m_j},
\eea
where $\bfm^{(J)}$ and  $A_{kn}^{(J)}$
are as in \eqref{A2def}.
With $m_{+}$ and $m_{++}$ as in \eqref{A2def}, let
\be\label{3M3}
\bfp_n^{(J)}= (p_{jn}^{(J)})_{J+1\le j\le n} 
=
\left(\frac{\Gamma(j-\alpha)/j!}
{\sum_{\ell=J+1}^n\Gamma(\ell-\alpha)/\ell!}\right)_{J+1\le j\le n},
\ee
and let 
 $\bfMult_{kn}^{(J)}$   
 be multinomial 
 $(k-  m_{+},\bfp_n^{(J)})$  with
\ben 
\PP\big(\bfMult_{kn}^{(J)}=(m_{J+1},\ldots, m_n) \big)=
(k- m_{+})!
\prod_{j=J+1}^n
\frac{(p_{jn}^{(J)})^{m_j}}{m_j!},
\een
for $m_j\ge 0$ and $\sum_{j=J+1}^n m_j=k-m_{+}$.
The summation over $\bfm^{(J)}$ of the LHS 
can be written as
\be\label{M9}
\PP\Big(\sum_{i=1}^{k-m_{+}} X_{in}^{(J)}=n- m_{++}\Big),
 \ee
 where $(X_{in}^{(J)})_{1\le i\le k-m_{+}}$ are i.i.d. with
\be\label{M8a}
\PP(X_{1n}^{(J)}=j) =p_{jn}^{(J)},\ J+1\le j\le n.
\ee
Again let $k'=k-m_{+}$ and $n'=n- m_{++}$ where convenient. 
Then 
\bea\label{3M01}
&&\sum_{\bfm^{(J)}\in A_{kn}^{(J)}}
\prod_{j=J+1}^n\frac{1}{m_j!} 
\Big(\frac{\Gamma(j-\alpha)}{j!}\Big)^{m_j} \cr
&&\cr
&&=
\frac{1}{k'!}
\sum_{\bfm^{(J)}\in A_{kn}^{(J)}}
k'!
\prod_{j=J+1}^n 
\frac{\big(p_{jn}^{(J)}\big)^{m_j}}{m_j!}
 \Big(\sum_{\ell=J+1}^n 
 \frac{\Gamma(\ell-\alpha)}{\ell!}\Big)^{m_j}\cr
 &&\cr
 &&=
\frac{1}{ k'!}
 \Big(\sum_{\ell=J+1}^n 
 \frac{\Gamma(\ell-\alpha)}{\ell!}\Big)^{k'}
\PP\Big(\sum_{i=1}^{k'} X_{in}^{(J)}=n'\Big).
\eea
To economise on notation we let here
\be\label{3defqL}
q_{jn}=\frac{\Gamma(j-\alpha)/j!}
{\sum_{\ell=1}^n\Gamma(\ell-\alpha)/\ell! }, \ j\in\N_J,
\ {\rm  so\ that}\ 
\frac{\sum_{\ell=J+1}^n \Gamma(\ell-\alpha)/\ell!}
{\sum_{\ell=1}^n\Gamma(\ell-\alpha)/\ell!}
=1- q_{+n},
\ee
still with $q_{+n}= \sum_{j=1}^J q_{jn}$.
(Previously, the $q_{jn}$ in \eqref{defqL} depended on $\lambda$, but those in \eqref{3defqL}  play an exactly analogous role.)
So we can write, from \eqref{M3.1}, 
\bea\label{3M4.2}
&&
\PP(M_{jn}(\alpha)=m_j, j\in\N_J,\, K_n(\alpha)=k)\cr
&&\cr
&&
\frac{1}{\alpha}
\frac{(k-1)!}{k'!  \prod_{j=1}^Jm_j!} 
\times
(1-q_{+n})^{k'}\prod_{j=1}^Jq_{jn}^{m_j}
\times
\Big(\frac{\alpha}{\Gamma(1-\alpha)} \sum_{\ell=1}^n 
 \frac{\Gamma(\ell-\alpha)}{\ell!}\Big)^k \cr
&&\cr
&&\hskip5cm
 \times
n\PP\Big(\sum_{i=1}^{k'} X_{in}^{(J)}=n'\Big).
\eea
Proceeding as in \eqref{M30a} we can write
\bea\label{3M30a}
&&
 \PP\Big(
\frac{M_{jn}(\alpha)}{K_n(\alpha)}\le 
q_j+a_j/n^{\alpha/2},\, j\in\N_J,\,
\frac{K_n(\alpha)}{n^\alpha} \le c\Big)\cr
&&\cr
&&= 
n^{\alpha(1+J/2)}  \int_{\bfy,x} x^J
\PP\big(\bfM_{jn}(\alpha)=
\lf (q_j+ y_j/n^{\alpha/2})xn^\alpha\rf,\, j\in\N_J, \cr
&&\cr
&& \hskip1.5cm  \times 
K_n(\alpha)=\lf xn^\alpha \rf\big)
\rmd \bfy\rmd x.
\eea

The following counterpart of Lemma \ref{lemm1} is proved  in the appendix.
\begin{lemma}\label{lemm2}
With the substitutions $k=\lf xn^\alpha \rf$, 
$m_j= (q_j+ y_j/n^{\alpha/2})k$, 
$k'=k-  m_{+}$,
we have the following limiting behaviour:
\bea\label{M1.32ab}
\frac{(k-1)!}{k'!\prod_{j=1}^Jm_j!}
(1-q_{+n})^{k'}\prod_{j=1}^Jq_{jn}^{m_j}
\sim
\frac{(xn^{\alpha})^{-J/2-1} e^{-\tfrac{x}{2} \bfy^T\bfQ^{-1}\bf y}}
{ \sqrt{(2\pi)^J{\rm det}(\bfQ)}};
\eea
\be\label{M1.44}
\lim_{n\to\infty}
\Big(\frac{\alpha}{\Gamma(1-\alpha)} \sum_{\ell=1}^n 
 \frac{\Gamma(\ell-\alpha)}{\ell!}\Big)^{\lf xn^\alpha\rf}
=
e^{-x/\Gamma(1-\alpha)};
\ee
\be\label{M1.43abc}
\lim_{n\to\infty}
n\PP\Big(\sum_{i=1}^{k'} X_{in}^{(J)}=n'\Big)
= f_{Y_x(1)}(1).
\ee
\end{lemma}

Lemma \ref{lemm2} provides the required components for \eqref{3M30a} and, keeping in mind the factor $n^{\alpha(1+J/2)} x^J$, \eqref{330} follows.
It remains to show how the density $f_{Y_x(1)}(1)$
in \eqref{330} is related to the density of the Mittag-Leffler distribution.
To do this we use a result of \cite{Covo2009a} which relates the L\'evy  density of a subordinator to that of a subordinator with the same but truncated L\'evy measure.
Write the probability density function (pdf) of a Stable($\alpha$) subordinator
$(S_x(\alpha))_{x\ge 0}$ (using variable $x\ge 0$ for the time parameter) having Laplace transform $e^{-x\tau^{\alpha}}$ and  L\'evy  density 
$\alpha z^{-\alpha-1}{\bf 1}_{\{z>0\}}/\Gamma(1-\alpha)$
as
\be\label{defSd}
f_{S_x(\alpha)}(s) =
\frac{1}{\pi} \sum_{j=0}^\infty 
\frac{(-1)^{k+1}}{k!}  
\frac{\Gamma(\alpha k+1)}{s^{\alpha k+1}}
x^k\sin(\pi\alpha k)
\ee
 (\cite{Pitman2006}, p.10).
Write the density of a  Mittag-Leffler random variable
 with parameter $\alpha$, $L_\alpha$, as
\be\label{defMLd}
f_{L_\alpha}(s) =
\frac{1}{\pi\alpha} \sum_{j=0}^\infty 
\frac{(-1)^{k+1}}{k!}  {\Gamma(\alpha k+1)}{s^{k-1}}
\sin(\pi\alpha k)
\ee
  (\cite{Pitman2006}, p.11),
and observe that
\be\label{SML}
\frac{1}{\alpha x}f_{S_x(\alpha)}(1)
=f_{L_\alpha}(x).
\ee
The L\'evy  density of $Y_x(1)$ is obtained from \eqref{LYP} with $\lambda=1$, and is the same as that of 
$S_x(\alpha)$ truncated at 1.
Using  Cor. 2.1 of \cite{Covo2009a} (in his formula set $x=1$, $s=\lambda$, $t=x$, and take 
$\Lambda(\lambda)=\lambda^{-\alpha}/\Gamma(1-\alpha)$) we have
\be\label{Kov}
f_{Y_x(\lambda)}(1)
=
e^{x\lambda^{-\alpha}/\Gamma(1-\alpha)}
\Big(f_{S_x(\alpha)}(1)
+\sum_{\kappa=1}^{\lc 1/\lambda \rc-1}
(-x)^\kappa A_{\lambda:\kappa}^{(1)}(1,x)\Big)
\ee
(with $\sum_1^0\equiv 0$), where the $ A_{\lambda:\kappa}^{(1)}$ are certain functions defined by Covo.
When $\lambda =1$ these functions disappear from the formula and we simply have
\ben
f_{Y_x(1)}(1)
=
e^{x/\Gamma(1-\alpha)}f_{S_x(\alpha)}(1).
\een
Using this together with \eqref{SML} to replace $f_{Y_x(1)}(1)$ in \eqref{330}
we obtain a representation of the limiting density of $K_n(\alpha)/n^\alpha$ in terms of a  Mittag-Leffler density, in agreement with Thm. 3.8 of \cite{Pitman2006}.\footnote{Pitman in fact gives almost sure convergence of  $K_n(\alpha)/n^\alpha$ to the  Mittag-Leffler,
whereas we have only convergence in distribution.
But we expand this to \eqref{330}.}

 \medskip\noindent {\bf Remark.}\
Efficient routines for numerical computation of one-sided L\'evy and Mittag-Leffler distributions are in \cite{SaaVen2011}.
   
 \section{The Ewens and Pitman Sampling Formulae}
 \label{s4}
 
  \subsection{The Ewens Sampling Formula}
 \label{ss44}
 In the next theorem, dealing with
  $\PD(\theta)$,   the limiting behaviours 
 of $\bfM_n(\theta)$ and $K_n(\theta)$ as Poisson and normal are well known  separately, but the asymptotic independence seems not to have been previously noted, and what's also interesting and potentially of further useful application is the way the joint limit arises from the methodology of the previous sections. 
 A result of  \cite{Hensley1982}  plays a key role.
 Let $m_j$ be nonnegative integers, $c\in \R$
  and write $(M_{jn}(\theta))_{1\le j\le n}$ for the components of $\bfM_n(\theta)$.
   
\begin{theorem}\label{M6}
For the $\PD(\theta)$ model we have the  asymptotic property:
\bea\label{4.330}
&&
\lim_{n\to\infty}
 \PP\Big(
M_{jn}(\theta)=m_j, j\in\N_J,
\frac{K_n(\theta)-\theta \log n}
{\sqrt{\theta\log n}} \le c\Big)
\cr
&&\cr
&& 
=
\prod_{j=1}^J
\frac{1}{m_j!} 
\Big(\frac{\theta}{j}\Big)^{m_j}e^{-\theta/j}
\times 
\frac{1}{\sqrt{2\pi }}
 \int_{x\le c} 
 e^{-\tfrac{1}{2} x^2}  \rmd x.
\eea
\end{theorem}
 
   \noindent {\bf Proof of Theorem \ref{M6}:}\
Starting from Ewens' sampling formula
\be\label{4.15}
\PP(\bfM_n(\theta)=\bfm,\, K_n(\theta)=k)
= \frac{n!\Gamma(\theta)\theta^k}{\Gamma(n+\theta)}
\prod_{j=1}^n 
\frac{1}{m_j!} 
\Big(\frac{1}{j}\Big)^{m_j}
\ee
(\cite{Ewens1972}, \cite{Pitman2006}, p.46), 
we get
\bea\label{5M3.1}
&&
\PP(M_{jn}(\theta)=m_j, j\in\N_J,\, K_n(\theta)=k)\cr
&&\cr
&&
=
 \frac{n!\Gamma(\theta)\theta^{k-  m_+}}{\Gamma(n+\theta)}
\prod_{j=1}^J\frac{1}{m_j!} 
\Big(\frac{\theta}{j}\Big)^{m_j}
\sum_{\bfm^{(J)}\in A_{kn}^{(J)}}
\prod_{j=J+1}^n
\frac{1}{m_j!} 
\Big(\frac{1}{j}\Big)^{m_j},
\eea
where $\bfm^{(J)}$ and  $A_{kn}^{(J)}$
are as in \eqref{A2def} and $m_+= \sum_{j=1}^J  m_j$.
Let
\be\label{5M3}
p_{jn} =\frac{1}{j}\Big/\sum_{\ell=J+1}^n\frac{1}{\ell},\
J+1\le j\le n,
\ee
and let  $\bfMult_{kn}:= (M_{jn})_{J+1\le j\le n}$ be multinomial 
with  distribution 
\ben 
\PP\big(\bfMult_{kn}=(m_{J+1},\ldots, m_n) \big)=
(k- m_+)!
\prod_{j=J+1}^n
\frac{(p_{jn})^{m_j}}{m_j!},
\een
for $m_j\ge 0$,  with
$\sum_{j=J+1}^n m_j=k- m_+$.
Thus, arguing as in \eqref{3M01}, we find
\be\label{5M9}
\sum_{\bfm^{(J)}\in A_{kn}^{(J)}}
\prod_{j=J+1}^n
\frac{1}{m_j!} 
\Big(\frac{1}{j}\Big)^{m_j}
=
\frac{1}{k'!}
\Big(\sum_{\ell=J+1}^n\frac{1}{\ell}\Big)^{k'}
\PP\Big(\sum_{i=1}^{k'} X_{in}=n'\Big),
 \ee
 where  $k'=k-  m_+$, $n'=n- m_{++}$
 and $(X_{in})_{1\le i\le k'}$ are i.i.d. with
\be\label{5M8a}
\PP(X_{1n}=j) =p_{jn},\ J+1\le j\le n.
\ee

So we can write, from \eqref{5M3.1}, 
\bea\label{4M4.2}
&&
\PP(M_{jn}(\theta)=m_j, j\in\N_J,\, K_n(\theta)=k)=\cr
&&\cr
&&=
 \frac{\Gamma(n)\Gamma(\theta)\theta^{k'}}{\Gamma(n+\theta)k'!}
 \Big(\sum_{\ell=J+1}^n\frac{1}{\ell}\Big)^{k'}
\prod_{j=1}^J\frac{1}{m_j!} 
\Big(\frac{\theta}{j}\Big)^{m_j}
n\PP\Big(\sum_{i=1}^{k'} X_{in}=n'\Big).
\eea
Then for $c\in \R$
\bea\label{430}
&&
 \PP\Big(
M_{jn}(\theta)=m_j, j\in\N_J,
\frac{K_n(\theta)-\theta \log n}
{\sqrt{\theta\log n}} \le c\Big)
\cr
&&\cr
&& 
=
 \int_{z\le \theta \log n +c\sqrt{\theta\log n}}
  \PP\big(M_{jn}(\theta)=m_j, j\in\N_J,
  K_n(\theta)=\lf  z\rf\big)
 \rmd z\cr
 &&\cr
 &&=\sqrt{\theta\log n}   \cr
 &&\cr
 &&\times
 \int_{x\le c}
  \PP\big(M_{jn}(\theta)=m_j,\, j\in\N_J,
  K_n(\theta)=\lf  \theta \log n +x\sqrt{\theta\log n}\rf\big)
 \rmd x.\cr
 &&
\eea
Let $k_n(x)= \lf  \theta \log n +x\sqrt{\theta\log n}\rf$ and calculate 
\bean
&&
 \sqrt{\theta\log n} \PP\big(M_{jn}(\theta)=m_j, j\in\N_J,
  K_n(\theta)=k_n(x)\big)\cr
  &&  \cr
  &&
  =
  \sqrt{\theta\log n}\times {\rm RHS\ of\  \eqref{4M4.2}\ with }\ k=k_n(x).
\eean
We let $n\to\infty$ in this. Consider the various factors.
First,
\be\label{4.31}
 \frac{\Gamma(n)\Gamma(\theta)}{\Gamma(n+\theta)}
\sim  \frac{\Gamma(\theta)}{n^\theta}.
\ee
Next, using a standard asymptotic for the harmonic series,
\bea\label{4.32}
 &&
  \Big(\sum_{\ell=J+1}^n\frac{1}{\ell}\Big)^{k- m_+}
  = \Big(\sum_{\ell=1}^n\frac{1}{\ell}
  -\sum_{j=1}^J \frac{1}{j}
 \Big)^{k-  m_+}\cr
 &&\cr
 &&
 =(\log n +\gamma- \sum_{j=1}^J 1/j+O(1/n)) ^{k-m_+}  \cr
 &&\cr
 &&
 =(\log n) ^{k- m_+}
 \Big(1+
 \frac{\gamma  -\sum_{j=1}^J 1/j+O(1/n)}{\log n}\Big)^{k- m_+}\cr
 &&\cr
 &&\sim
 (\log n) ^{k- m_+}
 e^{\theta\gamma-\theta\sum_{j=1}^J 1/j}, \ {\rm as}\ k=k_n(x)\to\infty.
\eea
Here $\gamma=0.566...$ is Euler's constant and we recall  $k_n(x)= \lf  \theta \log n +x\sqrt{\theta\log n}\rf $.

Substitute \eqref{4.32}  in \eqref{4M4.2}, remembering the factor of $\sqrt{\theta \log n}$ from \eqref{430},   to get  \eqref{4M4.2} asymptotic to
\bea\label{dubst}
\frac{ \sqrt{\theta\log n}}{n^\theta}
\frac{(\theta\log n)^{k'} }
{k'!}
\times
e^{\theta\gamma}
\Gamma(\theta)
\times
\prod_{j=1}^J
\frac{1}{m_j!} 
\Big(\frac{\theta}{j}\Big)^{m_j}
e^{-\theta/j}
\times
n\PP\Big(\sum_{i=1}^{k'} X_{in}=n'\Big).  \cr
&&
\eea
Using Stirlings's formula
and the  relations $k'=k-m_+=k- \sum_{j=1}^J  m_j$ and 
 $k=k_n(x)= \lf  \theta \log n +x\sqrt{\theta\log n}\rf $, we find
\bea\label{4.35}
&&\frac{ \sqrt{\theta\log n}}{n^\theta}
\frac{(\theta\log n)^{k'}}{k'!}
\sim
\frac{ \sqrt{\theta\log n}}{n^\theta}
\frac{(\theta\log n)^{k'}}
{\sqrt{2\pi k'}\, (k')^{k'}\, e^{-k'}
}\cr
 &&\cr
 &&\sim
 \frac{e^{k'}}{\sqrt{2\pi} n^\theta}
 \Big(\frac{\theta\log n}{k'}\Big)^{k'}
 =
 \frac{e^{k'}}{\sqrt{2\pi} n^\theta}
 \Big(1- \frac{k'-\theta\log n}{k'}\Big)^{k'}\cr
 &&\cr
 &&=
  \frac{e^{k'}}{\sqrt{2\pi} n^\theta}
 \exp\Big(k'\log\Big(1- \frac{k'-\theta\log n}{k'}\Big)\Big).
\eea
 Note that
 \ben
\frac{k'- \theta \log n}{k'}
= \frac{\lf  \theta \log n +x\sqrt{\theta\log n}\rf-  m_+-\theta\log n}{k'}
=
\frac{x\sqrt{\theta\log n}+O(1)}{k'}
\een
is $O\big(1/\sqrt{\log n}\big)\to 0$ as $n\to\infty$, 
so we can use the expansion  $\log (1-z)=-z-z^2/2-\cdots$ for small $z$ to get for the RHS of \eqref{4.35}
 \bean
 &&
  \frac{e^{k'}}{\sqrt{2\pi} n^\theta}
\exp\Big(k'\Big(-\frac{k'-\theta\log n}{k'}-
\frac{1}{2}
 \Big(\frac{k'-\theta\log n}{k'}\Big)^2\Big)\Big)
 +O\big(\frac{1}{\log n}\big)^{3/2}\cr
 &&\cr
 &&
 =
  \frac{1}{\sqrt{2\pi}}
  \exp\Big(-\frac{1}{2}
 \frac{ \big(x\sqrt{\theta\log n}+O(1)\big)^2}{ \theta \log n +x\sqrt{\theta\log n}} \Big)
 +O\big(\frac{1}{\log n}\big)^{1/2}\cr
 &&\cr
 &&
 \to
   \frac{1}{\sqrt{2\pi}} e^{-x^2/2}.
 \eean
Thus  \eqref{dubst} is asymptotic to
\be\label{4.40}
\prod_{j=1}^J
\frac{1}{m_j!} 
\Big(\frac{\theta}{j}\Big)^{m_j}
e^{-\theta/j}
\times
   \frac{1}{\sqrt{2\pi}} e^{-x^2/2}   
\times
   e^{\theta\gamma}
\Gamma(\theta)
\times
n\PP\Big(\sum_{i=1}^{k'} X_{in}=n'\Big).
 \ee

 It remains to deal with the $ X_{in}$ term  in \eqref{4.40}.
  Use the Fourier inversion formula for the mass function of  a discrete random variable (\cite{GK1968}, p.233) to write
   \bea\label{3.num1}
 n\PP\Big(\sum_{i=1}^{k'} X_{in}=n'\Big)
 &=&
 \frac{n}{2\pi}
  \int_{-\pi}^{\pi} e^{-n' \rmi \nu}
\big( \phi_{n}(\nu)\big) ^{  k'} \rmd \nu\cr
&&\cr
&=&
 \frac{n}{2\pi n'}
 \int_{-n'\pi}^{n'\pi} e^{-\rmi \nu}
\big( \EE\big(e^{\rmi\nu  X_{1n}/n'} \big) \big)^{  k'} \rmd \nu,
 \eea
 where $ \phi_{n}(\nu)= \EE\big(e^{\rmi\nu  X_{1n}})$, $\nu\in\R$.
By \eqref{5M8a}
\bean
 \EE\big(e^{\rmi\nu  X_{1n}/n'} \big)
 =
\sum_{j=J+1}^n \frac{e^{\rmi\nu j/n'}}{j}\Big/\sum_{\ell=J+1}^n\frac{1}{\ell}
 =1+
\sum_{j=J+1}^n \frac{1}{j}\big(e^{\rmi\nu j/n'}-1\big)
\Big/\sum_{\ell=J+1}^n\frac{1}{\ell},
\eean
in which 
\bean
&&\frac{\sum_{j=J+1}^n \big(e^{\rmi\nu j/n'}-1\big)/j}
{\sum_{\ell=J+1}^n/\ell}
=
\frac{\big(\sum_{j=1}^n-\sum_{j=1}^J\big)
 \big(e^{\rmi\nu j/n'}-1\big)/j
}
{\big(\sum_{\ell=1}^n-\sum_{\ell=1}^J \big)1/\ell}\cr
&&\cr
&&
=\frac{\sum_{j=1}^n \big(e^{\rmi\nu j/n'}-1\big)/j}
{\log n+O(1)}+ O\big(\frac{1}{n\log n}\big).\cr
\eean
The numerator here  is
\bean
&&
\sum_{j=1}^n \big(e^{\rmi\nu j/n'}-1\big)/j
=
\sum_{j=1}^n \frac{\rmi \nu}{n'}
\int_0^1
e^{\rmi\nu zj/n'}\rmd z
=
\frac{\rmi \nu}{n'}\int_0^1
\sum_{j=0}^n 
e^{\rmi\nu zj/n'}\rmd z-\frac{\rmi \nu}{n'}\cr
&&\cr
&&
=\frac{\rmi \nu}{n'}
\int_0^1 \frac{e^{\rmi\nu z(n+1)/n'}-1}
{e^{\rmi\nu z/n'}-1}\rmd z +O\big(\frac{1}{n}\big)
=\rmi \nu
\int_0^1 \frac{e^{\rmi\nu z(n+1)/n'}-1}
{\rmi\nu z+O(1/n)}\rmd z +O\big(\frac{1}{n}\big)\cr
&&\cr
&&
=
\int_0^1 \big(e^{\rmi\nu z}-1\big)\frac{\rmd z}{z}+O\big(\frac{1}{n}\big),
\eean
and consequently
\bean
&&
\Big(
1+
\sum_{j=J+1}^n \frac{1}{j}\big(e^{\rmi\nu j/n'}-1\big)
\Big/\sum_{\ell=J+1}^n\frac{1}{\ell}\Big)^{k_n(x)- m_+}\cr
&&\cr
&&
=
\Big(
1+
\frac{\int_0^1 \big(e^{\rmi\nu z}-1\big)\rmd z/z}
{\log n+O(1)}
\Big)^{\lf \theta \log n +x\sqrt{\theta\log n}- m_+\rf}.
\eean
The last expression has limit $\exp(\theta \int_0^1 (e^{\rmi\nu z}-1)\rmd z/z)$ so it follows that
\bean
\lim_{n\to\infty}
 \EE\big(e^{\rmi\nu  X_{1n}/n'} \big)^{k'}
 =e^{\theta \int_0^1 (e^{\rmi\nu z}-1)\rmd z/z}
 = \what w_\theta(\nu),
\eean
in the notation of \cite{Hensley1982}
(his Eqs. (2.6) and (2.12) with $\alpha=\theta$ and $K(\alpha)=1)$.
Hensley's function $\wh w_\theta(\nu)$ is the characteristic function of a density $w_\theta(u)$, so taking the limit  in \eqref{3.num1} gives
\ben
\frac{1}{2\pi}  \int_{-\infty}^{\infty} e^{-\rmi\nu}
\wh w_\theta(\nu)\rmd \nu =w_\theta(1).
\een
(Again justification for taking the limit through the integral in \eqref{3.num1} can be given, but is unnecessary.)
By \cite{Hensley1982}, Eq. (2.1), $w_\theta(1) = e^{-\theta\gamma}/\Gamma(\theta)$.
We conclude that
\ben 
\lim_{n\to\infty}
n\PP\Big(\sum_{i=1}^{k'} X_{in}=n'\Big)
= e^{-\theta\gamma}/\Gamma(\theta).
\een
Substituting this in \eqref{4.40}  gives \eqref{4.330}, clearly  a proper distribution.
\halmos

 \subsection{The Pitman Sampling Formula} \label{ss4}
 
The next theorem deals with Pitman's formula for  the distribution of  $\PD(\alpha, \theta)$. Introducing $\alpha$ changes the rate of growth of $K_n$ from logarithmic to power form, $n^\alpha$, as in $\PD_\alpha$.
Let $c\ge 0$, $J\ge 1$,
$\bfa=(a_1,a_2,\ldots,a_J)\in \R^J$, define $\bf Q$ as in \eqref{Qdef},
 and write $(M_{jn}(\alpha, \theta))_{1\le j\le n}$ for the components of $\bfM_n(\alpha,\theta)$.
 Recall that $f_{L_\alpha}$ is the Mittag-Leffler density
 in \eqref{defMLd}.
   
\begin{theorem}\label{M7}
For the $\PD(\alpha,\theta)$ model we have the  asymptotic property:
\bea\label{5.330}
&&
\lim_{n\to\infty}
 \PP\Big(
\frac{M_{jn}(\alpha, \theta)}{K_n(\alpha, \theta)}
\le 
q_j+ \frac{a_j}{n^{\alpha/2}}, j\in\N_J,
\frac{K_n(\alpha, \theta)}{n^\alpha} \le c\Big)
\cr
&&\cr
&& 
=
\frac{\Gamma(\theta+1)}{\Gamma(\theta/\alpha+1)}
  \int_{\bfy\le \bfa} \int_{x\le c} 
\frac{x^{J/2+\theta/\alpha}e^{-\tfrac{x}{2} \bfy^T\bfQ^{-1}\bf y}}
{\sqrt{(2\pi)^J{\rm det}(\bfQ)}}
f_{L_\alpha}(x)     \rmd x\rmd \bfy.\ \
\eea
\end{theorem}
    
\noindent {\bf Proof of Theorem \ref{M7}:}\
Most of the work for this is done in the
proof of Theorem \ref{M5}, so we just indicate the main  modifications needed. Start from Pitman's sampling formula
\begin{align}\label{FSF}
&
\PP(\bfM_n(\alpha, \theta)=\bfm,\, K_n(\alpha, \theta)=k)\cr
&\cr
&
= \frac{n!}{\alpha}
\Big(\frac{\alpha}{\Gamma(1-\alpha)}\Big)^k
\frac{\Gamma(\theta/\alpha+k)}
{\Gamma(\theta/\alpha+1)}
\frac{\Gamma(\theta+1)}{\Gamma(n+\theta)}
\times
\prod_{j=1}^n 
\frac{1}{m_j!} 
\Big(\frac{\Gamma(j-\alpha)}{j!}\Big)^{m_j}\ \
\end{align}
(\cite{PY1997}, Sect A.2, p.896, \cite{Pitman1995PTRF}, Prop. 9),
and proceed as in \eqref{3M30} to write
\bea\label{4M30}
&&
 \PP\Big(
\frac{M_{jn}(\alpha, \theta)}{K_n(\alpha,\theta)}\le 
q_j+a_j/n^{\alpha/2},\, j\in\N_J,
\frac{K_n(\alpha, \theta)}{n^\alpha} \le c\Big)\cr
&&\cr
&&= 
n^{\alpha(1+J/2)}  \int_{\bfy,x} x^J
\PP\big(\bfM_{jn}(\alpha, \theta)=
\lf (q_j+ y_j/n^{\alpha/2})xn^\alpha\rf,\,
 j\in\N_J, \cr
&&\cr
&& \hskip1.5cm  \times 
K_n(\alpha, \theta)=\lf xn^\alpha \rf\big)
\rmd\bfy \rmd x.
\eea
For the probability in \eqref{4M30} 
we get, as in \eqref{3M4.2}, 
\bea\label{4M4.2a}
&&
\PP(M_{jn}(\alpha, \theta)=m_j, \, j\in\N_J,\, K_n(\alpha, \theta)=k)
=\Big(\frac{\alpha}{\Gamma(1-\alpha)} \sum_{\ell=1}^n 
 \frac{\Gamma(\ell-\alpha)}{\ell!}\Big)^k \cr
&&\cr
&&
\times
 \frac{\Gamma(n)\Gamma(\theta/\alpha+k)\Gamma(\theta+1)}
{\alpha k'!\Gamma(\theta/\alpha+1)\Gamma(n+\theta)}
\big(1-q_{+n}\big)^{k'}
\prod_{j=1}^J \Big(\frac{q_{jn}^{m_j}}{m_j!}
\Big)
\times
 n\PP\Big(\sum_{i=1}^{k'} X_{in}^{(J)}=n'\Big),\cr
&&
\eea
where the $q_{jn}$ are defined in \eqref{3defqL} and
we replace $k$ by $\lf xn^\alpha\rf$ and
$m_j$ by  $(q_j+ y_j/n^{\alpha/2})k$.
Using $\Gamma(n)/\Gamma(n+\theta)\sim n^{-\theta}$
and $\Gamma(\theta/\alpha+k)\sim k^{\theta/\alpha}(k-1)!$ 
together with the estimates found in the
proof of Theorem \ref{M5}
for the other factors in \eqref{4M4.2a}, 
and recalling that \eqref{3M30} equals $f_{L_\alpha}(x)$,
we arrive at
\eqref{5.330}, in which the marginal limiting distribution for 
$K_n(\alpha,\theta)/n^\alpha$ agrees with that in Eq. (3.27), p.68,  of \cite{Pitman2006}.
 \halmos

 \section{Other Limits}\label{s5} 
In this section we derive other limits of the processes studied above.
Subsection \ref{LB} analyses $\PD_\alpha^{(r)}$ as $r\dto 0$ or $r\to\infty$.
Subsection \ref{LC} summarises the results in  a  convergence diagram for $\bfM_n(\alpha,r)$,  showing that the convergences as $n\to\infty$ and $r \dto 0$ commute.

 \subsection{Limiting behaviour of  $\PD_\alpha^{(r)}$ as $r\dto 0$ or $r\to\infty$}\label{LB}
 Kingman's   $\PD_\alpha$ distribution arises by taking the ordered jumps, $(\Delta S_1^{(i)})_{i\ge 1}$, up till time 1, of a driftless $\alpha$-stable subordinator $S=(S_t)_{t>0}$, normalised by their sum $S_1$,  as  a random distribution on the infinite simplex
  $\nabla_{\infty}$.
A natural generalisation  
is to delete the $r$ largest jumps ($r\ge 1$ an integer) up till time 1
 of   
$S$  and consider  the vector whose components are  the remaining jumps $(\Delta S_1^{(i)})_{i\ge r+1}$ normalised by their sum.
Again we obtain  a random distribution on $\nabla_{\infty}$, now with an extra parameter, $r$.
This is  the $\PD_\alpha^{(r)}$ distribution of  \cite{IpsenMaller2017}. Some limiting and other properties of it were studied in  \cite{IpsenMallerShemehsavar2020},
 where it was extended to all $r>0$, not just integer $r$,
 and in \cite{CZ2023}.

  By its construction 
$\PD_\alpha^{(r)}$ reduces to  $\PD_\alpha$ when $r$ is set equal to $0$, but Theorem \ref{1} shows a kind of continuity property:
as $r\dto 0$, $r>0$, the distribution of $(\bfM_n(\alpha,r), K_n(\alpha,r))$
from  $\PD_\alpha^{(r)}$ converges to the corresponding quantity from  $\PD_\alpha$; a useful result for statistical analysis, where we might fit the  $\PD_\alpha^{(r)}$ model with general $\alpha$ and $r$ to data, and test $H_0:r=0$, i.e., whether $r$ is needed in the model to describe the data.
The parameter $r$ models  ``over-dispersion'' in the data 
(\cite{CZ2023}).

 In this subsection we keep $n$ and $\alpha$ fixed and let $r\dto 0$ or $r\to\infty$.
 Part (i) of the next theorem is an analogue of the  \cite{PY1997}, p.880, result that, for each $\theta>0$, the
limit of $\PD(\alpha,\theta)$ as $\alpha\dto 0$ is 
$\PD(\theta)$.
 
 \begin{theorem}\label{1}
 We have the following limiting behaviour for  $\PD_\alpha^{(r)}$.
 
 Part  (i) \ As $r\dto 0$, the limit of $\PP(\bfM_n(\alpha,r)=\bfm,\, K_n(\alpha,r)=k)$ is 
 the distribution in \eqref{3.15}, i.e., the  distribution of $(\bfM_n(\alpha), K_n(\alpha))$.
 
 
 
 Part  (ii) \ In the limit as $r\to \infty$,  a sample  of
  $\PD_\alpha^{(r)}$
is of $n$ species each with 1 representative,
 i.e., the singleton distribution with mass $ {\bf 1}_{(\bfm,k)=((n,0,\ldots,0_n), n)}$.
 \end{theorem}
 

\noindent {\bf Proof of  Theorem \ref{1}:}
For this proof it's convenient to
 work from a formula
 for $\PD_\alpha^{(r)}$ given by 
  \cite{IpsenMallerShemehsavar2020}:
\bea\label{7}
&&
\PP(\bfM_n(\alpha,r)=\bfm,\, K_n(\alpha,r)=k) \cr
&&\cr
&&
=  
n \int_{0}^{\infty}\ell_n(\lambda)^k
 \PP\Big({\rm \bf Mult}(k,\bfp_n(\lambda))= \bfm\Big)
\PP\Big( {\rm Negbin}\big(r, \frac{1}{\Psi(\lambda)}\big) =k\Big)
\frac{\rmd \lambda}{\lambda}. \cr &&
\eea
Here
$ \bfm\in A_{kn}$ (see \eqref{4.4a0})
and $k  \in  \N_n$. The function $\ell_n(\lambda)$ is defined  as 
\be\label{OL6}
\ell_n(\lambda) 
:=
\frac{\Psi_n(\lambda)-1}{\Psi(\lambda)-1}
=
\frac{ \int_0^\lambda\sum_{j=1}^n (z^j/j!) z^{-\alpha-1}e^{-z}\rmd z}
{\int_0^\lambda z^{-\alpha-1}(1-e^{-z})\rmd z}
\le 1,
\ee
where $\Psi(\lambda)$  is as in \eqref{2.7aa}, 
${\rm \bf Mult}(k,\bfp_n(\lambda))$ is multinomial $(k,\bfp_n(\lambda))$  with
\be\label{OL3}
\bfp_n(\lambda)= (p_{jn}(\lambda))_{j\in\N_n} =
\left(\frac{F_j(\lambda)}{\sum_{\ell=1}^n F_\ell(\lambda)}\right)_{j\in\N_n},
\ee
and $F_j(\lambda)$ is as in \eqref{Fdef}.
Finally, in \eqref{7}, 
$ {\rm Negbin}\big(r, 1/\Psi(\lambda)\big)$ is a negative binomial rv with parameter $r>0$ and success probability $1/\Psi(\lambda)$,
thus 
\be\label{OL22}
\PP\Big( {\rm Negbin}\big(r, \frac{1}{\Psi(\lambda )}\big) =
k  \Big)
= 
\frac{\Gamma(r+ k)}{\Gamma(r) 
k!}
\left( 1- \frac{1}{\Psi(\lambda )}\right)^{k}
\Big(\frac{1}{\Psi(\lambda) }\Big)^r,\ k\in \N_n.
\ee

\noindent
{\it Part  (i)}\ 
    Observe, by \eqref{2.7aa}, for all $\lambda>0$, 
 \be\label{5a}
1\vee \lambda^\alpha 
\int_{0}^{\lambda}e^{-z}z^{-\alpha} \rmd z
\le
\Psi(\lambda)=1+\alpha\lambda^\alpha \int_{0}^\lambda
(1-e^{-z})z^{-\alpha-1} \rmd z
\le 1+\lambda^\alpha\Gamma(1-\alpha),
\ee
so
\be\label{25}
\PP\Big( {\rm Negbin}\big(r, \frac{1}{\Psi(\lambda )}\big) =k\Big)
=
\frac{\Gamma(r+ k)}{\Gamma(r) k!}
\frac{\big(\Psi(\lambda )-1\big)^k }
{\Psi(\lambda )^{k+r}}
\le
\frac{\Gamma(r+ k)\Gamma^k(1-\alpha)}{\Gamma(r) k!} \lambda^{k\alpha}.
\ee

Fix $A>0$.  Using \eqref{25} and $\ell_n(\lambda)\le 1$ the 
component of the 
integral over $(0,A)$ in \eqref{7} does not exceed
 \ben
 n\int_0^A  \frac{\Gamma(r+ k)\Gamma^k(1-\alpha)}{\Gamma(r) k!}
\lambda^{k\alpha-1}  \rmd \lambda
 =
 \frac{n\Gamma(r+ k)\Gamma^k(1-\alpha)}{\Gamma(r) k!}
 \times \frac{A^{k\alpha} }{k\alpha},
\een
 and this tends to 0 as $r\dto 0$ since $\Gamma(r)\to\infty$ then.

For the component of the integral over $(A,\infty)$ in \eqref{7}, change variable to $\lambda=\lambda/r$ and write it  as
 \be\label{10} 
n \int_{Ar}^{\infty}\ell_n(\lambda/r)^k
  \PP\Big({\rm \bf Mult}(k,\bfp_n(\lambda/r))= \bfm\Big)
\PP\Big( {\rm Negbin}\big(r, \frac{1}{\Psi(\lambda/r)}\big) =k\Big)
 \, \frac{\rmd \lambda}{\lambda}. 
\ee
Define quantities $ \ell_n^\infty$ and $\bfM_{kn}^\infty$ so that
$ \ell_n^\infty=\lim_{\lambda\to\infty} \ell_n(\lambda)$
and 
\ben
\PP\Big(\bfM_{kn}^\infty= \bfm\Big)=
\lim_{\lambda\to\infty}  \PP\Big({\rm \bf Mult}(k,\bfp_n(\lambda))= \bfm\Big).
\een 
Given $\veps>0$, choose $A=A(\veps)$ large enough for these limits to be within $\veps$ of $\ell_n(\lambda)$ and the multinomial probability when $\lambda>A$,  
and note that, by \eqref{5a},
 \ben 
\Psi(\lambda/r)
\ge
(\lambda/r)^\alpha 
\int_{0}^{\lambda/r}e^{-z}z^{-\alpha} \rmd z
\ge
(\lambda/r)^\alpha 
\int_{0}^{A}e^{-z}z^{-\alpha} \rmd z
=: c_{\alpha,A}(\lambda/r)^\alpha
\een
when $\lambda/r>A$. 
So an upper bound for the integral in \eqref{10} is, for $r<1$,
 \bea\label{26} 
&& \frac{n} {c_{\alpha,A}^r} \int_{Ar}^{\infty}\left(
(\ell_n^\infty)^k  \PP\big(\bfM_{kn}^\infty= \bfm\big)+\veps\right)
\frac{\Gamma(r+ k)}{\Gamma(r) k!} 
\lambda^{-\alpha r -1}r^{\alpha r }
\rmd \lambda\cr
 &&\cr
 &&=
  \frac{n} {c_{\alpha,A}^r}\left(
(\ell_n^\infty)^k  \PP\big(\bfM_{kn}^\infty= \bfm\big)+\veps\right)
\frac{\Gamma(r+ k)}{\alpha r\Gamma(r) k!} 
A^{-\alpha r}.
\eea
Let $r\dto 0$ and note that
$c_{\alpha,A}^r\to 1$, $A^{-\alpha r}\to 1$, 
$ r\Gamma(r) =\Gamma(r+1) \to 1$ and $\Gamma(r+ k)\to \Gamma( k)=(k-1)!$.
Also from \eqref{OL6}, \eqref{OL3} and \eqref{Fdef},
\ben
\ell_n^\infty =\frac{\alpha}{\Gamma(1-\alpha)} 
\sum_{j=1}^n \frac{\Gamma(j-\alpha)} {j!}
\een
and 
\ben
p_j(\infty) =\left(\frac{F_j(\infty)}{\sum_{\ell=1}^n F_\ell(\infty)}\right)
= \frac{\Gamma(j-\alpha)} {j!}\slash 
\sum_{\ell=1}^n \frac{\Gamma(\ell-\alpha)} {\ell!},
\een
so that
\ben
(\ell_n^\infty)^k
  \PP\big(\bfM_{kn}^\infty=\bfm, K_n(\infty)=k) =
  k!     \Big(\frac{\alpha}{\Gamma(1-\alpha)}\Big)^k
\prod_{j=1}^n \frac{1}{m_j!} \Big(\frac{\Gamma(j-\alpha)}{j!}\Big)^{m_j}.
\een
Finally letting $\veps\dto 0$  we obtain an upper bound for the limit of the integral in \eqref{7} 
equal to  the distribution in \eqref{3.15}.

For a lower bound, we can discard the 
component of the 
integral over $(0,A)$ in \eqref{7}, then change variables as in \eqref{10} and use Fatou's lemma.
Thus we prove Part (i).

 \smallskip\noindent
 {\it Part  (ii)}\  In  the integral  in \eqref{7}, change variable to 
 \ben
 \rho = \big(1- \frac{1}{\Psi(\lambda)}\big)
 \frac{r}{1/\Psi(\lambda)} = r(\Psi(\lambda)-1),
 \een
 so that 
 \ben
\frac{\rmd \lambda}{\lambda}
=\frac{\rmd \rho}{r\lambda\Psi'(\lambda)}, 
\quad
 \Psi(\lambda)= \frac{\rho}{r}+1\
\ {\rm and}\ \
 \lambda =\Psi^{\leftarrow}\Big( \frac{\rho}{r}+1\Big),
 \een
 where $\Psi^{\leftarrow}$ is the inverse function to $\Psi$.
 Then $r\to\infty$ implies $\Psi(\lambda)\to 1$ thus $\lambda\to 0$.
 From \eqref{OL6} and \eqref{OL3} and L'Hopital's rule we see that
\ben
\lim_{\lambda\to 0} 
\ell_n(\lambda)=1
\  {\rm and}\
\lim_{\lambda\to 0} 
p_{jn}(\lambda) = p_{jn}(0)=
\lim_{\lambda\to 0} 
\frac{1}{j!\sum_{\ell=1}^n  \lambda^{\ell-j}/\ell!}
 = {\bf 1}_{j=1}.
\een
After  changing variable in \eqref{7} we get
\bea\label{72}
&&
\PP(\bfM_n(\alpha,r)=\bfm,\, K_n(\alpha,r)=k) \cr
&&\cr
&&
=  
n \int_{0}^{\infty}\ell_n\Big(\Psi^{\leftarrow}\big( \frac{\rho}{r}+1\big)\Big)^k
 \PP\Big(
 {\rm \bf Mult}(k,\bfp_n (\Psi^{\leftarrow}\big( \frac{\rho}{r}+1\big))= \bfm\Big) \cr
 &&\cr
 &&\hskip3cm
  \times 
\PP\Big( {\rm Negbin}\big(r, \frac{r}{r+\rho}\big) =k\Big)
 \, \frac{\rmd \rho}{r\lambda \Psi'(\lambda)  },
\eea
in which, by \eqref{2.7aa},
$\Psi(\lambda)-1 \sim \alpha \lambda/(1-\alpha)$ as $\lambda\to 0$, so 
$r\lambda \sim (1-\alpha)\rho/\alpha$, and 
 $\Psi'(\lambda)=
  \alpha \int_{0}^1e^{-\lambda z}z^{-\alpha} \rmd z
 \to \alpha/(1-\alpha)$ as $\lambda\to 0$.
As $r\to \infty$ we have convergence of the negative binomial to the Poisson, so the RHS of \eqref{72} converges to 
\bean
&&
n \int_{0}^{\infty}
 \PP\Big(
 {\rm \bf Mult}(k,\bfp_n (0))= (1,0,\ldots 0)\Big) 
 \frac{e^{-\rho} \rho^k}{k!}
  \frac{ \rmd \rho}{\rho }\cr
  &&\cr
  &&
  =\frac{n}{k}  \PP\big({\rm \bf Mult}(k,\bfp_n(0))= (1,0,\ldots 0)\big).
\eean
Here ${\rm \bf Mult}(k,\bfp_n(0)$ is multinomial $(k, \bfp_n(0)$ with
$\bfp_n(0)= (p_{jn}(0)_{1\le j\le n}$, 
and so the probability equals $ {\bf 1}_{(\bfm,k)=((n,0,\ldots,0_n), n)}$ as required.
\halmos

\begin{theorem}\label{5.2}
For the  $\PD_\alpha^{(r)}$  model we have the  property:
\ben 
\lim_{r\dto 0} ({\rm RHS\ of}\ \eqref{30})
=
{\rm RHS\ of}\ \eqref{330}.
\een
\end{theorem}

\noindent {\bf Proof of  Theorem \ref{5.2}:}
From \eqref{30} it suffices to consider the integrand
\be\label{sci}
 \frac{x^{r+J/2-1} }{\Gamma(r)\Gamma^r(1-\alpha)}
 \frac{e^{-\tfrac{x}{2} \bfy^T\bfQ^{-1}\bf y}}
{\sqrt{(2\pi )^J {\rm det}(\bfQ)}}
\int_{\lambda=0}^\infty
e^{-x(\lambda^{-\alpha}\vee 1)/\Gamma(1-\alpha)}
f_{Y_x(\lambda)}(1)
\frac{\rmd \lambda}{ \lambda^{\alpha r+1}},
\ee
and for this we look at $\int_0^1$ and $\int_1^\infty$ separately.
Recall that $f_{Y_x(\lambda)}(y)$ is uniformly bounded in $\lambda$  by $C$, say.
For the first part, make the transformation $y=\lambda^{-\alpha}$, 
so  $\lambda=y^{-1/\alpha}$ and 
$\rm d\lambda=-y^{-1/\alpha-1}\rmd y/\alpha$.  Then  
$\int_0^1$ is bounded by
\bean
&&
 \frac{Cx^{r+J/2-1} }{\Gamma(r)\Gamma^r(1-\alpha)}
 \frac{e^{-\tfrac{x}{2} \bfy^T\bfQ^{-1}\bf y}}
{\sqrt{(2\pi )^J {\rm det}(\bfQ)}}
\int_{\lambda=0}^1
e^{-x(\lambda^{-\alpha}\vee 1)/\Gamma(1-\alpha)}
\frac{\rmd \lambda}{ \lambda^{\alpha r+1}} \cr
&&\cr
&&
=
 \frac{Cx^{r+J/2-1} }{\alpha\Gamma(r)\Gamma^r(1-\alpha)}
 \frac{e^{-\tfrac{x}{2} \bfy^T\bfQ^{-1}\bf y}}
{\sqrt{(2\pi )^J {\rm det}(\bfQ)}}
\int_{y=1}^\infty
e^{-xy/\Gamma(1-\alpha)}
y^{r-1} \rmd y,
\eean
which  tends to 0 as $r\dto 0$. So we can neglect $\int_0^1$.
For  $\int_1^\infty$ the contribution is
\bean
&&
 \frac{x^{r+J/2-1} }{\Gamma(r)\Gamma^r(1-\alpha)}
 \frac{e^{-\tfrac{x}{2} \bfy^T\bfQ^{-1}\bf y}}
{\sqrt{(2\pi )^J {\rm det}(\bfQ)}}
\int_{\lambda=1}^\infty
e^{-x/\Gamma(1-\alpha)}
f_{Y_x(1)}(1)
\frac{\rmd \lambda}{ \lambda^{\alpha r+1}} \cr
&&\cr
&&
=
 \frac{x^{r+J/2-1} }
 {\alpha r\Gamma(r)\Gamma^r(1-\alpha)}
 \frac{e^{-\tfrac{x}{2} \bfy^T\bfQ^{-1}\bf y}}
{\sqrt{(2\pi )^J {\rm det}(\bfQ)}}
 e^{-x/\Gamma(1-\alpha)}
f_{Y_x(1)}(1)\cr
&&
\to
\frac{x^{J/2-1} e^{-\tfrac{x}{2} \bfy^T\bfQ^{-1}\bf y}}
{\alpha\sqrt{(2\pi )^J {\rm det}(\bfQ)}}
 e^{-x/\Gamma(1-\alpha)}
f_{Y_x(1)}(1),
\eean
using that $ r\Gamma(r)=\Gamma(r+1)\to 1$ as $r\dto 0$.
This equals the integrand on the  RHS of \eqref{330}.   \halmos

  \medskip\noindent {\bf Remarks.}\  
  (i)\ We can simplify the limit distribution in \eqref{30} somewhat by using 
  Cor. 2.1 of \cite{Covo2009a} 
  as we did in the  proof of Theorem \ref{M5}.
  Once again in the integral in \eqref{sci}  look at $\int_0^1$ and $\int_1^\infty$ separately.
  The component of the integral over $(1,\infty)$ can be explicitly evaluated and equals 
  \ben
  \frac{1}{\Gamma(r)\Gamma^r(1-\alpha)} x^r 
  f_{L_\alpha}(x).
  \een
The component over $(0,1)$ involves the functions $ A_{\lambda:\kappa}^{(1)}$ defined by Covo (see \eqref{Kov}),  which can be calculated using a recursion formula he gives; or, for an alternative approach, see \cite{perman1993}.
  
  (ii)\  Theorem 3.1 of \cite{MallerShemehsavar2023} shows that a sampling model based on the Dickman subordinator has the same limiting behaviour as the Ewens sampling formula.
  This is no surprise as the  Dickman function is well known to be closely associated with $\PD(\theta)$ and some of its generalisations;
  see for example
\cite{ABT2003}, pp.14, 76, 
  \cite{AB2015},
  \cite{Handa2009}, 
\cite{IpsenMallerShemehsavar2020,IpsenMallerShemehsavar2021},
  \cite{MallerShemehsavar2023}.

   \subsection{Convergence Diagram for   $\PD_\alpha^{(r)}$}\label{LC}
  Using Theorems \ref{M4} and \ref{M5}
  we can give the following
   convergence diagram:

\begin{figure}[H]
\centering
\begin{tikzpicture}
\node  (ul)  at (0,4)  {$(\bfM_n(\alpha,r), K_n(\alpha,r))$  };
\node  (ur)  at (5,4)  {$(\bfM(\alpha,r), K(\alpha,r))$};
\node  (ll)  at (0,1)  { $(\bfM_n(\alpha), K_n(\alpha))$};
\node  (lr)  at (5,1)  {$(\bfM(\alpha), K(\alpha))$};

\draw[->] 
(ul) edge node [midway, sloped, above] {$n\to \infty$} 
          node [midway, sloped, below] {}(ur)
(ul) edge node [midway, sloped, above] {$r\to 0$}
      node [midway, sloped, below]  {} (ll)
(ur) edge node [midway, sloped, above] {$r\to 0$} 
          node [midway,  below] {\,\,\,\,\,\,\,}(lr)
(ll) edge node [midway, sloped, above] {$n\to \infty$} 
          node [midway, sloped, below] {}(lr)   ;
\end{tikzpicture}
\caption{{\it Convergence diagram for   $\PD_\alpha^{(r)}$.
Upper left to lower left
and 
lower left to lower right are in Theorems  4.1 and 3.1 of \cite{MallerShemehsavar2023}.
For upper left to upper right,  $(\bfM(\alpha,r), K(\alpha,r))$ has the distribution on the RHS of  \eqref{30}.
For upper right to lower right, the convergence is proved in Theorem \ref{M5}, with $(\bfM(\alpha),K(\alpha))$ denoting a random variable with the  Normal-Mittag-Leffler  mixture  distribution on the RHS of  \eqref{330}.\\
The diagram is schematic only: the random variables have to be normed and centered in certain  ways before convergence takes place.
}
}\label{diagramM}
\end{figure}

We can also use Theorem \ref{M4} 
to complete the 
   convergence diagram for $\bfM_n(\alpha,r)$
   in \cite{MallerShemehsavar2023} in which $r\to\infty$, $\alpha\dto 0$, such that $r\alpha\to a>0$; the missing  upper right entry in their Figure 2  is  $\bfM(\alpha,r)$.
%
%
  We also have the following  useful results which follow immediately from  Theorems \ref{M4} and \ref{M5}.

  \begin{corollary}\label{cor1}
     $\sqrt{ K(\alpha,r)}\bfM(\alpha,r) \eqdr N(\bf0,\bfQ)$, independent of  $K(\alpha,r)$; and 
   $\sqrt{ K(\alpha)}\bfM(\alpha) \eqdr N(\bf0,\bfQ)$, independent of  $K(\alpha)$.
  \end{corollary}

%
 
 \section{Discussion}\label{s6}
The methods developed in Sections \ref{s2}--\ref{s5} offer a unified approach  with nice interpretations to limit theorems of various sorts for Poisson-Dirichlet models, and we expect 
can be applied in other situations too. 
 For example, \cite{GSpeed2002} analyse a ``general selection model" based on  the infinite alleles mutation scheme.
 Their formula (1.13)
  gives a version of the distribution of $(\bfM_n, K_n)$ amenable to analysis by our methods.
\cite{RWF2013} (see their Eqs. (3.1), (3.2))
analyse a species sampling model based on  normalized  inverse-Gaussian diffusions as a generalisation of $\PD_\alpha$;  an extra parameter $\beta$ is introduced to the $\PD_\alpha$ model, somehat analogous to our $\PD_\alpha^{(r)}$.
 See also \cite{LijoiMenaPrunster2005} and \cite{FF2014}.

Formula \eqref{4.1} can be compared with an analogous formula, Eq. (4.14), p.81, of \cite{Pitman2006}, in which the component 
$\prod_{i=1}^k \Psi^{(n_i)}(\lambda)$
can be replaced by $\prod_{j=1}^n \big(\Psi^{(j)}(\lambda)\big)^{m_j}$ using
 the identity
$m_j =\sum_{i=1}^k 1_{\{j=n_i\}}$.
See also \cite{Pitman1997}
and the Gibbs partitions in \cite{Pitman2006}, pp.25-26, Eq. (1.52).
\cite{Hansen1994} gave a general treatment of decomposable combinatorial structures having a Poisson-Dirichlet limit.

There are many other limiting results in the literature.
\cite{JoyceKroneKurtz2002} prove a variety of
Gaussian limit theorems 
as the mutation rate $\theta\to \infty$ in the $\PD(\theta)$ model;
see also \cite{Griffiths1979}, \cite{Handa2009}.
 \cite{Feng2007,Feng2010} gives large deviation results as $\theta\to\infty$.
  \cite{LancelotJames2008} gives results like this related to $\PD(\alpha,\theta)$  in a Bayesian setup.
\cite{DF2020} investigate $\alpha$-diversity in $\PD(\alpha,\theta)$. 
  When $r  \to  \infty$, size-biased samples from 
  $\PD_\alpha^{(r)}$
tend to those of $\PD(1 -  \alpha )$, i.e., to $\PD(\theta)$ with $\theta = 1- \alpha $.
Together with Part (i) of Theorem \ref{1} this suggests that  $\PD_\alpha^{(r)}$ is intermediate in a sense between $\PD_\alpha$ and $\PD(\theta)$.
\cite{CZ2023} show that  $\PD_\alpha^{(r)}$ can alternatively be obtained as the ordered normalised jumps of 
a negative binomial process (\cite{Gregoire1984}).
When $r  \to  \infty$,  
$\alpha   \to  0$ such that $r \alpha   \to  a > 0$, a limit involving the Dickman subordinator
is obtained from   $\PD_\alpha^{(r)}$
(\cite{MallerShemehsavar2023}).

The appearance of the normal as  the conditional  limiting distribution of $\bfM_n$ given  $K_n$, after  centering and norming, for 
  $\PD_\alpha^{(r)}$, $\PD_\alpha$ and $\PD(\alpha,\theta)$,
is useful for statistical applications.
In \cite{IpsenMallerShemehsavar2018} we used a least squares method to fit
 $\PD_\alpha^{(r)}$ to a variety 
 of data sets, some quite small. Even so we observed by simulations that the finite sample distributions of the estimates of $\alpha$ and $r$ were quite closely approximated by the normal.
 The asymptotic normality derived in  the present paper provides a rationale for using least squares and helps explain the approximate normality of
 the parameter estimates. Similar ideas may be useful for 
  $\PD_\alpha$ and $\PD(\alpha,\theta)$.
Functionals of $\bfM_n$ 
such as $\sum_{j=1}^n (M_{jn}- EM_{jn})^2/EM_{jn}$
are important in many aspects of population genetics; see for example the discussion of measures of sample diversity  in Section 4 of
\cite{Watterson1974}.

 \section{Appendix: Proofs of Lemmas \ref{lemm1} and \ref{lemm2}}\label{app}

 \subsection{Proof of Lemma \ref{lemm1}}\
  To prove \eqref{L1}, consider
 the first factor in  \eqref{M4.2}, in which we let 
 $y_{jn}= y_j/n^{\alpha/2}$, 
 $h_{jn} = q_j+ y_{jn}$ and 
$h_{+n} =\sum_{j=1}^J h_{jn}$; 
then
replace $m_j$ by  $h_{jn}k$, 
recall $k'=k-m_{+}
=(1-h_{+n})k$,
and calculate
\bea\label{M1.3}
&&\frac{\Gamma(r+k)}{k'!\prod_{j=1}^Jm_j!}\sim
\frac{
\sqrt{2\pi (r+k-1)}\, (r+k-1)^{r+k-1} e^{-(r+k-1)} 
}{\sqrt{2\pi(1- h_{+n})k}
\,  ( (1-h_{+n})k)^{ (1-h_{+n})k} e^{- (1-h_{+n})k}
}\cr
&&\cr
&&\hskip4cm
 \times
\frac{1}
{
\prod_{j=1}^J
\sqrt{2\pi h_{jn}k}\,
 (h_{jn}k)^{h_{jn}k} e^{-h_{jn}k}
}, 
\eea
using Stirling's formula for the factorials.
Some calculations reduce this to
\be\label{M1.31}
\frac{
k^{r-J/2-1}}
{(2\pi)^{J/2} \sqrt{ (1- h_{+n})\prod_{j=1}^J h_{jn} }
}    \times
\frac{1}{ (1-h_{+n})^{ (1-h_{+n})k} 
\prod_{j=1}^J h_{jn}^{h_{jn}k} 
}.  
\ee
Multiply this by the $q$-terms  in \eqref{M4.2}
 to get
\bea\label{M1.32}
&&\frac{\Gamma(r+k)}{k'!\prod_{j=1}^Jm_j!}
(1-q_{+n})^{k- m_{+}}\prod_{j=1}^Jq_{jn}^{m_j}
  \cr
&&\cr
&&
\sim
\frac{
k^{r-J/2-1}}
{\sqrt{(2\pi)^{J} (1- h_{+n})\prod_{j=1}^J h_{jn}}
}     \times
\Big(\frac{1-q_{+n}}
{1-h_{+n}}\Big)^{ (1-h_{+n})k} \prod_{j=1}^J 
\Big(\frac{q_{jn}}{ h_{jn}}\Big)^{h_{jn}k},\cr
&&
\eea
where  $q_{+n}= \sum_{j=1}^J q_{jn}$.
Note 
from \eqref{defqL} that, as $n\to\infty$, 
\be\label{fqL}
q_{jn}\to 
\frac{F_j(\infty)} {\sum_{\ell=1}^\infty F_\ell(\infty)}
=\frac{\Gamma(j-\alpha)/j!}
{\sum_{\ell=1}^\infty\Gamma(\ell-\alpha)/\ell! }
=\frac{\alpha\Gamma(j-\alpha)}
{j!\Gamma(1-\alpha) }=q_j, \ j\in\N_J,
\ee
 and
 $1-q_{+n}\to 1-q_+$, where $q_+ =\sum_{j=1}^Jq_j$.
Also 
 $h_{jn} = q_j+ y_{jn} \to q_j$.
 Thus
 \bea\label{M1.33}
 &&
\Big(\frac{1-q_{+n}}
{1-h_{+n}}\Big)^{ (1-h_{+n})k} \times
\prod_{j=1}^J 
\Big(\frac{q_{jn}}{ h_{jn}}\Big)^{h_{jn}k}
\cr
&&\cr
&& =\Big(1+ \frac{\sum_{j=1}^J (y_{jn}+q_j-q_{jn})}
{
1-\sum_{j=1}^J (q_j+y_{jn}) }\Big)
^{(1-\sum_{j=1}^J  (q_j+y_{jn}))k}    \cr
&&\cr
&&
\times
\prod_{j=1}^J 
\Big(1- \frac{
y_{jn}+q_j-q_{jn}
}{
q_j+y_{jn}
}\Big)
^{(q_j+y_{jn})k}\cr
&&\cr
&&=
\exp\Big(
\sum_{j=1}^J (q_j+y_{jn}) k
\log\Big(1- \frac{
y_{jn}+q_j-q_{jn}
}
{q_j+y_{jn}}
\Big)\cr
&&\cr
&&
+\Big(1-\sum_{j=1}^J (q_j+y_{jn})\Big)
k\log\Big(1+ \frac{
\sum_{j=1}^J (y_{jn}+q_j-q_{jn})
}
{
1-\sum_{j=1}^J (q_j+y_{jn})
}
\Big).
  \eea
 Expand the RHS of \eqref{M1.33}
using $\log (1-z)=-z-z^2/2-\cdots$
 and   $\log (1+z)=z-z^2/2-\cdots$ for small $z$.
 The first order terms cancel and we're left with
  \bea\label{M1.34}
 &&
- \frac{k}{2}\Big(
\frac{\sum_{j=1}^J 
(y_{jn}+q_j-q_{jn})^2}{q_j+y_{jn}}
+\frac{\big(\sum_{j=1}^J (y_{jn}+q_j-q_{jn}\big)^2}
{1-\sum_{j=1}^J (q_j+y_{jn})}
\Big).
 \eea
 Recalling $k=\lf xn^\alpha\rf$,  $y_{jn}=y_j/n^{\alpha/2}$ and $q_j-q_{jn}=O(n^{-\alpha})$, the last expression simplifies to be asymptotic to
   \be\label{M1.35}
- \frac{1}{2}xn^\alpha
 \Big(
 \frac{1}{n^\alpha}
 \sum_{j=1}^J 
\frac{y_j^2}{ q_j} 
+ 
\frac{(\sum_{j=1}^J y_j)^2}
{n^\alpha (1-q_+)}\Big) \to
- \frac{x}{2} 
 \Big(
 \sum_{j=1}^J 
\frac{y_j^2}{ q_j} 
+ 
\frac{(\sum_{j=1}^J y_j)^2}
{1-q_+}\Big).
 \ee
The RHS can be written in the form
 $-x \bfy^T \bfQ^{-1}\bf y/2$,
 so the limit of the LHS of  \eqref{M1.33}
 is  $\exp(-x \bfy^T \bfQ^{-1}\bf y/2)$,
 and thus \eqref{M1.32} is asymptotic to
 \be\label{M1.48}
 \frac{k^{r-J/2-1}e^{-\tfrac{x}{2} \bfy^T\bfQ^{-1}\bf y}}
{ \sqrt{(2\pi)^{J}
(1-\sum_{j=1}^Jq_j)\prod_{j=1}^Jq_j}}
=
\frac{k^{r-J/2-1}e^{-\tfrac{x}{2} \bfy^T\bfQ^{-1}\bf y}}
{\sqrt{(2\pi)^{J}{\rm det}(\bfQ)}}.
\ee
(See  Lemma \ref{lemm.determinant} in the Appendix for calculating  the determinant of $\bfQ$.)

To prove \eqref{L2}, 
consider the second factor in  \eqref{M4.2}.
From \eqref{2.7aa} and $\int_0^\infty (1-e^{- z})\alpha z^{-\alpha-1} \rmd z=\Gamma(1-\alpha)$ we see that
\bea\label{M1.36}
&&
\Psi(\lambda n)=1+ (\lambda n)^\alpha
\int_{0}^{\lambda n} (1-e^{- z})\alpha z^{-\alpha-1} \rmd z\cr
&&=
1+ (\lambda n)^\alpha\Gamma(1-\alpha)- (\lambda n)^\alpha
\int_{\lambda n}^\infty (1-e^{- z})\alpha z^{-\alpha-1} \rmd z,
\eea
in which, integrating by parts,
\ben
 (\lambda n)^{\alpha}
\int_{\lambda n}^\infty (1-e^{- z})\alpha z^{-\alpha-1} \rmd z
= 1-e^{-\lambda n}
 + (\lambda n)^{\alpha}\int_{\lambda n}^\infty e^{- z} z^{-\alpha} \rmd z
 = 1+O(e^{-\lambda n}).
 \een
So   $ \Psi(\lambda n) = (\lambda n)^{\alpha}\Gamma(1-\alpha)
+o(n^{-\alpha})$, and consequently
\bea\label{3M1.37}
\frac{(\lambda n)^\alpha}{\Psi(\lambda n)}
=\frac{1}{\Gamma(1-\alpha)} +o(n^{-\alpha}).
\eea
To estimate $\sum_{\ell=1}^n F_\ell(\lambda n)$, recall \eqref{Fdef} and write it as 
\bea\label{M1.38}
\sum_{\ell=1}^n F_\ell(\lambda n) 
=
\alpha\sum_{\ell=1}^n \frac{\Gamma(\ell-\alpha)}{\ell!} 
-\alpha \sum_{\ell=1}^n \int_{\lambda n}^\infty
\frac{z^\ell}{\ell!} z^{-\alpha-1}e^{-z} \rmd z.
\eea
The first term on the RHS of \eqref{M1.38} is 
\bean
&&\alpha\Big(\sum_{\ell=1}^\infty - \sum_{\ell>n}\Big) 
\int_{0}^\infty
\frac{z^\ell}{\ell!} z^{-\alpha-1}e^{-z} \rmd z\cr
&&=
\alpha \int_{0}^\infty
(1-e^{-z}) z^{-\alpha-1} \rmd z
- \alpha \int_{0}^\infty \sum_{\ell>n}\frac{z^\ell}{\ell!} z^{-\alpha-1}e^{-z} \rmd z
\eean
The first term on the RHS here equals $\Gamma(1-\alpha)$.
The second term equals
\ben
n^{-\alpha} \int_{0}^\infty
\sum_{\ell>n}\frac{(nz)^\ell}{\ell!} 
\alpha z^{-\alpha-1}e^{-nz} \rmd z=
n^{-\alpha} \int_{0}^\infty
\PP\Big(\frac{{\rm Poiss}(nz)}{nz}>\frac{1}{z}\Big)\alpha z^{-\alpha-1} \rmd z.
\een
Since $\PP\big({\rm Poiss}(nz)>n\big) \le nz/n=z$ we can apply dominated convergence and  the weak law of large numbers to see that
\be\label{3M1.35}
\lim_{n\to\infty}
\int_{0}^1
\PP\Big(\frac{{\rm Poiss}(nz)}{nz}>\frac{1}{z}\Big)\alpha z^{-\alpha-1} \rmd z
= 0, 
\ee
and similarly
\ben
\lim_{n\to\infty}
 \int_1^\infty
\PP\Big(\frac{{\rm Poiss}(nz)}{nz}>\frac{1}{z}\Big)\alpha z^{-\alpha-1} \rmd z
=
 \int_1^\infty \alpha z^{-\alpha-1} \rmd z
=1.
\een 
 It follows that 
\be\label{3M1.37b}
\alpha\sum_{\ell>n}  \frac{\Gamma(\ell-\alpha)}{\ell!} 
=n^{-\alpha}+ o(n^{-\alpha}),
\ee
and consequently
\be\label{3M1.36b}
\alpha\sum_{\ell=1}^n \frac{\Gamma(\ell-\alpha)}{\ell!} 
=\Gamma(1-\alpha) -n^{-\alpha}+ o(n^{-\alpha}).
\ee
Returning to \eqref{M1.38}, we next estimate
\bean
&&
\alpha \sum_{\ell=1}^n \int_{\lambda n}^\infty
\frac{z^\ell}{\ell!} z^{-\alpha-1}e^{-z} \rmd z
=
\alpha\Big(\sum_{\ell=1}^\infty - \sum_{\ell>n}\Big) 
\int_{\lambda n}^\infty
\frac{z^\ell}{\ell!} z^{-\alpha-1}e^{-z} \rmd z\cr
&&\cr
&&=\alpha 
n^{-\alpha} \int_{\lambda}^\infty
(1-e^{-nz}) z^{-\alpha-1} \rmd z\cr
&& \hskip4cm 
- n^{-\alpha}
\int_{\lambda}^\infty
\PP\Big(\frac{{\rm Poiss}(nz)}{nz}>\frac{1}{z}\Big)
\alpha z^{-\alpha-1}\rmd z.
\eean
The first term on the RHS here equals 
$\lambda^{-\alpha} n^{-\alpha} +o(n^{-\alpha})$.
In the second term, whenever $\lambda<1$ the component of the integral over
$[\lambda,1]$  is $o(n^{-\alpha})$  by \eqref{3M1.35}, and the component  over $[1,\infty]$ is $1+o(1)$, and when $\lambda>1$, the integral equals $\lambda^{-\alpha}+o(1)$.
Thus  this second term is 
$\lambda^{-\alpha}n^{-\alpha}{\bf 1}_{\{\lambda> 1\}}
+ n^{-\alpha}{\bf 1}_{\{\lambda\le 1\}}+o(n^{-\alpha})$.
Subtracting this from the first term we obtain
\bean
&&
\alpha \sum_{\ell=1}^n \int_{\lambda n}^\infty
\frac{z^\ell}{\ell!} z^{-\alpha-1}e^{-z} \rmd z
=\lambda^{-\alpha}n^{-\alpha}{\bf 1}_{\{\lambda\le 1\}}
- n^{-\alpha}{\bf 1}_{\{\lambda\le 1\}}+o(n^{-\alpha}),
\eean
and subtracting this from the RHS of  \eqref{M1.38}, keeping in mind \eqref{3M1.36b}, we get
\bea\label{M1.38a}
\sum_{\ell=1}^n F_\ell(\lambda n) 
&=&
\Gamma(1-\alpha) -n^{-\alpha}-
\lambda^{-\alpha}n^{-\alpha}{\bf 1}_{\{\lambda\le 1\}}
+ n^{-\alpha}{\bf 1}_{\{\lambda\le 1\}}+o(n^{-\alpha})\cr
&&\cr
&=&
 \Gamma(1-\alpha) -n^{-\alpha}\big(\lambda^{-\alpha}\vee 1\big) +o(n^{-\alpha}).
\eea
Combining this with \eqref{3M1.37} and recalling $k=\lf xn^\alpha \rf$, we obtain
\bea\label{M1.40}
\Big(\frac{(\lambda n)^\alpha}{\Psi(\lambda n)}
\sum_{\ell=1}^n F_\ell(\lambda n)\Big)^{k}
&=&
\Big(1- \frac{\big(\lambda^{-\alpha}\vee 1\big)/
\Gamma(1-\alpha)}{n^{\alpha}}  +o(n^{-\alpha})\Big)^{\lf xn^\alpha \rf}\cr
&&\cr
&\to&
e^{-x (\lambda^{-\alpha}\vee 1)/
\Gamma(1-\alpha)},\ {\rm as}\ n\to \infty.
\eea

 To prove \eqref{L3}: the third factor in \eqref{M4.2} is 
 $1/\Gamma(r)\Psi(\lambda n)^{r}$, and by \eqref{3M1.37}
\be\label{nefa}
\frac{1}{\Gamma(r)\Psi(\lambda n)^{r}}
\sim
\frac{1}{\Gamma(r)(\lambda n)^{\alpha r}\Gamma^r(1-\alpha)}.
\ee

This just leaves the fourth factor in \eqref{M4.2}.
To prove \eqref{L4} we modify the local limit result obtained in  IMS, where we dealt with a sequence $( X_{in})$
only slightly different to $ (X_{in}^{(J)})$; namely, 
rather then \eqref{M8}, we had
$\PP(X_{1n}(\lambda)=j) =p_{jn}(\lambda)$,
$1\le j\le n$, where, rather then \eqref{M3}, we had
$\bfp_n(\lambda)= (p_{jn}(\lambda))_{1\le j\le n}$,
with the summation on the RHS of \eqref{M3} being over $1\le \ell \le n$ rather than $J+1\le \ell\le n$. We indicate just the main modifications needed.

First, we claim Prop. 3.1 of  IMS remains true if
$X_{in}(\lambda n)$ therein is replaced by $X_{in}^{(J)}(\lambda n)$, the denominator $\Gamma(1-\alpha)$ in Eq. (3.10) of  IMS
 is replaced by 
 $\Gamma(1-\alpha)\big(1-q_+\big)$, 
and corresponding modifications are made   in Eq. (3.11) and 
Eq. (3.12) of  IMS.
To check this, note that $p_{jn}(\lambda)$ in the proof of 
 Prop. 3.1 of  IMS 
 is replaced by $p_{jn}^{(J)}(\lambda)$, so 
 in Eq. (3.10) of  IMS, once $hn>3$, we need only replace the denominator by
 \bean  
&&\sum_{j=J+1}^n F_j(\lambda n) 
=
\sum_{j=1}^n F_j(\lambda n)
-\sum_{j=1}^{J} F_j(\lambda n)
\cr
&=&
\Gamma(1-\alpha)-\sum_{j=1}^{J} \Gamma(j-\alpha)/j!
+o(1)
=\Gamma(1-\alpha)\big(1-q_+\big)+o(1).
\eean
Here, note that $ F_j(\lambda n)\to \alpha \Gamma(j-\alpha)/j!=
\Gamma(1-\alpha)q_j$ by \eqref{Fdef} and \eqref{qdef}, apply  \eqref{M1.38a} and recall that $q_+= \sum_{j=1}^{J}q_j$.
The rest of the proof  of Part (a) of  Prop. 3.1 of  IMS  remains the same so the only modification necessary is to replace the denominator 
$\Gamma(1-\alpha)$ by  $\Gamma(1-\alpha)(1-q_+)$.
Likewise  the proofs  of Part (b) and (c) of  the proposition remain valid after corresponding modifications.

We used these results in  Prop. 3.2 of  IMS to show that the limiting distribution of 
$n^{-1} \sum_{i=1}^{\lfloor xn^\alpha\rfloor} X_{in}(\lambda n)$
is an 
 infinitely divisible distribution with characteristic exponent
 $\pi_\lambda(\rmd y)$ given by  Eq. (2.7) of  IMS.
 For the present situation,  $\pi_\lambda(\rmd y)$ is modified to
  $\pi_\lambda^{(J)} (\rmd y)$ just by replacing
   the denominator 
$\Gamma(1-\alpha)$  in Eq. (2.7) of  IMS with  $\Gamma(1-\alpha)(1-q_+)$.
Thus, letting $\phi_{\lambda n}^{(J)}(\nu) 
:= \EE\exp\big(\rmi\nu  X_{1n}^{(J)}(\lambda n)\big)$, $\nu\in\R$, and recalling that   $k'=k-m_{+} \sim
 \lf xn^\alpha (1-q_+)\rf$, 
 we have 
%
 \bean 
 \lim_{n\to\infty}
  \big(\phi_{\lambda n}^{(J)}(\nu/n') \big)^{k'}
&=&
\lim_{n\to\infty} 
 \EE\exp\Big(\frac{ \rmi\nu}{n'} \sum_{i=1}
 ^{ \lf xn^\alpha (1-q_+)\rf} X_{in}^{(J)}(\lambda n)\Big)\cr
 &&\cr
&=&
 e^{-  x(1-q_+)\int_{\R\setminus\{0\}}
 (e^{ \rmi\nu y} -1 ) \Pi_\lambda^{(J)}(\rmd y)}\cr
 &&\cr
 & =&
  e^{-  x  \int_{\R\setminus\{0\}}
 (e^{ \rmi\nu y} -1 ) \Pi_\lambda(\rmd y)}
 =\EE (e^{ \rmi\nu  Y_x(\lambda)}).
\eean
Note that the factor of  $(1-q_+)$
cancels and \eqref{L4}
follows as in ISM.
This completes the proof of Lemma \ref{lemm1}. \halmos

 \subsection{Proof of Lemma \ref{lemm2}}\
To find the limit as $n\to\infty$ of the probability in \eqref{3M30}
 we use  \eqref{3M4.2}, in which as in the proof of Lemma \ref{lemm1}
 we   let  $y_{jn}= y_j/n^{\alpha/2}$, 
 $h_{jn} = q_j+ y_{jn}$,
replace $m_j$ by  $h_{jn}k$, 
and consider first the  factor
\bean 
&&\frac{(k-1)!}{k'!\prod_{j=1}^Jm_j!}
=\frac{(k-1)!}
{ ((1-h_{+n})k)!
\prod_{j=1}^J (h_{jn}k)!}\cr
&&\cr
&&
\sim
\frac{
k^{-J/2-1}}
{ \sqrt{(2\pi)^{J} (1- h_{+n})\prod_{j=1}^J h_{jn} }
}    \times
\frac{1}{ (1-h_{+n})^{ (1-h_{+n})k} 
\prod_{j=1}^J h_{jn}^{h_{jn}k} 
}. 
\eean
(see \eqref{M1.3}).
Multiply this by the $q$-terms  in \eqref{3M4.2}
 to get
\bean 
&&\frac{(k-1)!}{k'!\prod_{j=1}^Jm_j!}
(1-q_{+n})^{k'}\prod_{j=1}^Jq_{jn}^{m_j}
  \cr
&&\cr
&&
\sim
\frac{
k^{-J/2-1}}
{\sqrt{(2\pi)^{J} (1- h_{+n})\prod_{j=1}^J h_{jn}}
}    
\Big(\frac{1-q_{+n}}
{1-h_{+n}}\Big)^{ (1-h_{+n})k} 
\prod_{j=1}^J 
\Big(\frac{q_{jn}}{ h_{jn}}\Big)^{h_{jn}k}.
\eean
Very similar working as in \eqref{M1.33}--\eqref{M1.48} then gives the asymptotic in 
\eqref{M1.32ab}   
 for this, 
and similar working as in \eqref{M1.36} 
 and \eqref{3M1.36b} proves  \eqref{M1.44}.
 
  Finally in \eqref{3M4.2}  is the factor involving the $ X_{in}$. 
  Using  Fourier inversion as in \eqref{3.num1}
 \be\label{3.num0}
 n\PP\Big(\sum_{i=1}^{k'} X_{in}^{(J)}=n'\Big)
=
 \frac{n}{2\pi n'} \int_{-n'\pi}^{n'\pi} e^{-\rmi \nu}
\big( \phi_{n}^{(J)}(\nu/n')\big) ^{  k'} \rmd \nu,
 \ee
where  $\phi_{n}^{(J)}(\nu) 
:= \EE\big(\exp(\rmi\nu  X_{1n}^{(J)})\big)$, $\nu\in\R$.
For this we have
\bea\label{724}
&&\phi_{n}^{(J)}(\nu/n') = 
\frac{
\sum_{j=J+1}^n e^{\rmi \nu j/n' } \Gamma(j-\alpha)/j!
}
{\sum_{\ell=J+1}^n \Gamma(\ell-\alpha)/\ell!
}\cr
&&\cr
&&=
\frac{
\sum_{j=1}^n e^{\rmi \nu j/n' } \Gamma(j-\alpha)/j!
- \sum_{j=1}^J e^{\rmi \nu j/n' } \Gamma(j-\alpha)/j!
}
{\sum_{\ell=1}^n \Gamma(\ell-\alpha)/\ell!
-\sum_{\ell=1}^J \Gamma(\ell-\alpha)/\ell!
}.
\eea
In  the numerator of \eqref{724}  we can replace 
the exponentials in the second summation 
 by $1+O(1/n)$, and we will next show that
\be\label{3.num}
\sum_{j=1}^n \frac{ e^{\rmi \nu j/n} \Gamma(j-\alpha)}{j!}
=
\frac{\Gamma(1-\alpha)}{\alpha} +
n^{-\alpha} \int_0^1(e^{\rmi \nu z}-1)z^{-\alpha-1}\rmd z
-\frac{n^{-\alpha}}{\alpha} +o(n^{-\alpha}). 
\ee
(For notational  simplicity, in \eqref{3.num} and what follows, we
replace $n'$ by $n$, which is irrelevant asymptotically.)
To prove \eqref{3.num},  first consider
\bea\label{4.num}
&&
\sum_{j=1}^\infty
 \frac{ e^{\rmi \nu j/n } \Gamma(j-\alpha)}{j!}
=
\int_0^\infty
\sum_{j=1}^\infty \frac{ e^{\rmi \nu j/n } z^j} {j!}
e^{-z} z^{-\alpha-1}\rmd z\cr
&&=
\int_0^\infty
 \big(e^{ze^{\rmi \nu /n}}-1\big)e^{-z} 
 z^{-\alpha-1}\rmd z
 \cr
&&=
\int_0^\infty
 \big(e^z-1\big)e^{-z}  z^{-\alpha-1}\rmd z
 +
 \int_0^\infty
 \big(e^{ze^{\rmi \nu /n}}-e^z\big)e^{-z} 
 z^{-\alpha-1}\rmd z\cr
&&\cr
&&
=\int_0^\infty
 \big(1-e^{-z}) z^{-\alpha-1}\rmd z
 +n^{-\alpha}  \int_0^\infty
 \big(e^{nz(e^{\rmi \nu /n}-1)}-1\big) 
 z^{-\alpha-1}\rmd z\cr
 &&\cr
 &&=
 \Gamma(1-\alpha)/\alpha 
 + n^{-\alpha} \int_0^\infty(e^{\rmi \nu z}-1)z^{-\alpha-1}\rmd z
  +o(n^{-\alpha}).
\eea
Next, towards \eqref{3.num}, we show
\be\label{5.num}
\sum_{j>n}
 \frac{ e^{\rmi \nu j/n } \Gamma(j-\alpha)}{j!}
= \frac{n^{-\alpha}}{\alpha}
 + n^{-\alpha} \int_1^\infty(e^{\rmi \nu z}-1)z^{-\alpha-1}\rmd z
 +o(n^{-\alpha}).
\ee
Write the LHS of \eqref{5.num} as
\bea\label{6.num}
\sum_{j>n} \frac{\Gamma(j-\alpha)}{j!}
+
\sum_{j>n}
\big( e^{\rmi \nu j/n }-1\big) \frac{ \Gamma(j-\alpha)}{j!}.
\eea
The first term in \eqref{6.num} equals
$\Gamma(1-\alpha) -n^{-\alpha}+ o(n^{-\alpha})$, by  \eqref{3M1.36b}. 
To deal with the second term, define rvs $Y_n$ with
\be\label{7.num}
P(Y_n=j)=
\frac{\Gamma(j-\alpha)}{
j!\sum_{\ell>n}  \Gamma(\ell-\alpha)/\ell!},\  j>n,
\ee
and note that $Y_n>n$ w.p.1. Take $u>1$ and calculate
\bea\label{7.num}
&&
P(Y_n/n\le u)= P(Y_n\le nu)= P(n<Y_n\le nu)\cr
&&\cr
&&
=
\sum_{n<j\le nu}\frac{\Gamma(j-\alpha)}{j!}
\Big/\sum_{j>n}\frac{\Gamma(j-\alpha)}{j!}\cr
&&\cr
&&
=
\alpha n^{\alpha} (1+o(1))
\sum_{n<j\le nu}
\frac{\Gamma(j-\alpha)}{j!},
\eea
where we used \eqref{3M1.37b} to estimate the denominator.
Then note that
\bea\label{8.num}
&&\alpha n^{\alpha}\sum_{n<j\le nu}\frac{\Gamma(j-\alpha)}{j!}
=
\alpha\int_0^\infty\sum_{n<j\le nu}
\frac{(nz)^j}{j!} e^{-nz} z^{-\alpha-1}\rmd z\cr
&&\cr
&&
=
\alpha\int_0^\infty
P\big(n<{\rm Poiss}(nz)\le nu\big)z^{-\alpha-1}\rmd z.
\eea
This expression equals
\ben
\int_0^\infty
P\Big(
\frac{n(1-z)}{\sqrt{nz}}
<\frac{{\rm Poiss}(nz)-nz}{\sqrt{nz}}
\le \frac{n(u-z)}{\sqrt{nz}}\Big)\alpha
z^{-\alpha-1}\rmd z,
\een
and when $0<z<1$, $P({\rm Poiss}(nz)>n)\le z$, so we can apply dominated convergence to see that the last expression  converges to
\bea\label{9.num}
\int_0^\infty
P\big(-\infty<N(0,1)<\infty)
{\bf 1}_{\{1<z\le u\}}\alpha
z^{-\alpha-1}\rmd z
=1-u^{-\alpha}.
\eea
Consequently
$Y_n/n \todr Y$, where
$P(Y\le u) = (1-u^{-\alpha}) {\bf 1}_{\{u>1\}}$ and so
\be\label{10.num}
E \big(e^{\rmi \nu Y_n/n }-1\big)
\to
\int_1^\infty  e^{\rmi \nu u }\alpha u^{-\alpha-1}\rmd u-1
=
\int_1^\infty\big(e^{\rmi \nu u }-1\big)\alpha u^{-\alpha-1}\rmd u.
\ee
It follows that 
\bea\label{11.num}
&&
\sum_{j>n}
\big( e^{\rmi \nu j/n }-1\big) \frac{ \Gamma(j-\alpha)}{j!}
\sim
\frac{n^{-\alpha}}{\alpha}\Big(1+ E \big(e^{\rmi \nu Y_n/n }-1\big)
\Big)\cr
&&\cr
&&
=
\frac{n^{-\alpha}}{\alpha} + n^{-\alpha}
\int_1^\infty\big(e^{\rmi \nu z }-1\big)\alpha z^{-\alpha-1}\rmd z
+o(n^{-\alpha}), 
\eea
proving \eqref{5.num}.
Subtracting \eqref{5.num}  from \eqref{4.num} gives \eqref{3.num}.

The denominator in \eqref{724} is,  using \eqref{3M1.36b},
\bea\label{3.numa}
&&
\Gamma(1-\alpha)/\alpha 
-\sum_{\ell=1}^J \Gamma(\ell-\alpha)/\ell!
-n^{-\alpha}/\alpha
+o(n^{-\alpha})\cr
&&=
\Gamma(1-\alpha)(1-q_+)/\alpha-n^{-\alpha}/\alpha+o(n^{-\alpha}),
\eea
and dividing this into \eqref{3.num},  the ratio of numerator and denominator is
\bea\label{3.numb}
1+
\frac{ \int_0^1(e^{\rmi \nu z}-1)\alpha z^{-\alpha-1}\rmd z
/\Gamma(1-\alpha)}
{n^\alpha(1- q_+)-1}  +o(n^{-\alpha}).
\eea
Raised to power $k'=k-m_{+}\sim xn^\alpha(1- q_+)$, this   converges to 
\ben
x \int_0^1(e^{\rmi \nu z}-1)\alpha z^{-\alpha-1}\rmd z
/\Gamma(1-\alpha).
\een
Note that again the factor  $1-q_+$ cancels,
 and it follows as in ISM that 
\be\label{M1.43a}
\lim_{n\to\infty}
n\PP\Big(\sum_{i=1}^{k'} X_{in}^{(J)}=n'\Big)
= f_{Y_x(1)}(1).
\ee

This completes the proof of Lemma \ref{lemm2}. \halmos

\begin{lemma}\label{lemm.determinant}
The  matrix $\bfQ$ defined via \eqref{Qdef} has
determinant 
\ben
(1-\sum_{j=1}^Jq_j)\prod_{j=1}^Jq_j.
\een
\end{lemma}

 \medskip\noindent{\bf Proof of Lemma \ref{lemm.determinant}}\
 Denote by $\bfA$ the $K\times K$ matrix in  \eqref{Qdef}
 (without the premultiplying factor), let 
 $\bfD= {\rm(diag}(Q_1-Q, \ldots, Q_K-Q)$,  let $\bf1$ be a 
 $K$-vector of 1s, and let ${\rm adj} (\bfD)$ be the  adjugate matrix of $\bfD$.
 Then
  $\bfA=\bfD+Q\bf1\bf1^T$ and by the Sherman-Morrison formula
 \ben
  {\rm det}(\bfA) 
=
   {\rm det} (\bfD) +Q\bf1^T {\rm adj}(\bfD){\bf 1}
   =
{\rm   \prod_{j=1}^J (Q_j-Q) +Q\sum_{i=1}^J \prod_{j\ne i}(Q_j-Q).}
   \een
 Let $a= 1-\sum_{j=1}^J q_j$ and note that
 \ben
 Q_i = (1-\sum_{j\ne i}q_j)\prod_{j\ne i}q_j
 =\frac{(a+q_i)Q}{q_i},
 \een
 so $Q_i-Q = aQ/q_i$. 
 Hence
 \bean
 && \prod_{j=1}^J (Q_j-Q) +Q\sum_{i=1}^J \prod_{j\ne i}(Q_j-Q)
=
   \prod_{j=1}^J (Q_j-Q)\Big( 1+Q \sum_{i=1}^J\frac{q_i}{aQ}\Big)\cr 
   &=&
      \prod_{j=1}^J (Q_j-Q)\Big( 1+Q\frac{1-a}{aQ}\Big)
= \frac{1}{a} \prod_{j=1}^J \frac{aQ}{q_j}
      \cr 
      &=&
      \frac{a^JQ^J}{aQ}= a^{J-1}Q^{J-1}.
 \eean
Taking into account  the premultiplying factor  in  \eqref{Qdef}
gives  $a^{-1}Q^{-1}$ for the determinant of the matrix in   \eqref{Qdef} and hence $aQ = ( 1-\sum_{j=1}^J q_j) \prod_{j=1}^Jq_j$ as the determinant of $Q$.   \halmos

{}


\begin{thebibliography}{}
%



\bibitem[\protect\citeauthoryear{Arratia, Barbour \& Tavar\'e}{Arratia, Barbour \& Tavar\'e}{2003}]{ABT2003}
Arratia, R., Barbour, A. \& Tavar\'e, S. (2003)
\newblock {\em Logarithmic Combinatorial Structures: A Probabilistic Approach,}
 EMS Monographs in Mathematics, Euro. Math. Soc., Zurich.

\bibitem[\protect\citeauthoryear{Arratia \& Baxendale}{Arratia \& Baxendale}{2015}]{AB2015}
Arratia, R. \& Baxendale, P.  (2015)
\newblock Bounded size bias coupling: a Gamma function bound,
and universal Dickman-function behavior.
\newblock {\em Probab. Theor.  Rel. Fields},  162, 411--429.


%





%




\bibitem[\protect\citeauthoryear{Chegini \& Zarepour}{Chegini \& Zarepour}{2023}]{CZ2023}
Chegini,  S. \& Zarepour, M. (2023)
\newblock Random discrete probability measures based on negative binomial process.
\newblock {\em arxiv.org/pdf/2307.00176}.

\bibitem[\protect\citeauthoryear{Covo}{Covo}{2009a}]{Covo2009a}
Covo, S.  (2009a)
\newblock One-dimensional distributions of subordinators with upper truncated Lévy measure, and applications.
\newblock {\em Adv. Appl. Prob.}, 41, 367–392.

\bibitem[\protect\citeauthoryear{Covo}{Covo}{2009b}]{Covo2009b}
Covo, S.  (2009b)
\newblock On approximations of small jumps of subordinators with particular emphasis on a Dickman-type limit.
\newblock {\em J. Appl. Prob.}, 46, 732--755.

%



%



\bibitem[\protect\citeauthoryear{Dolera \& Favaro}
{Dolera \& Favaro}{2020}]{DF2020}
Dolera, E. \& Favaro, S. (2020)
\newblock {A Berry-Eesseen theorem for Pitman’s $\alpha$-diversity}
\newblock  {\it Ann. Appl. Prob.}, 30, 847--869.

\bibitem[\protect\citeauthoryear{Ewens}{Ewens}{1972}]{Ewens1972}
Ewens, W. (1972)
\newblock The sampling theory of selectively neutral alleles.
\newblock {\em Theoret. Pop. Biol.},  3, 87--112.

\bibitem[\protect\citeauthoryear{Favaro \& Feng}{Favaro \& Feng}{2014}]{FF2014}
Favaro, S. \& Feng, S. (2014) 
\newblock Asymptotics for the number of blocks in a conditional Ewens-Pitman sampling model
\newblock {\it Electron. J. Probab. 19 (2014), no. 21, 1–15.}


\bibitem[\protect\citeauthoryear{Feng}{Feng}{2007}]{Feng2007}
Feng, S. (2007)
\newblock Large deviations associated with Poisson-Dirichlet
distribution and Ewens sampling formula.
\newblock {\it Ann.  Appl. Prob.}, 17, 1570--1595.


\bibitem[\protect\citeauthoryear{Feng}{Feng}{2010}]{Feng2010}
Feng, S. (2010)
\newblock {\em The {Poisson-Dirichlet} Distribution and Related Topics: Models and Asymptotic Behaviours}.
\newblock Probability and its Applications. Springer.

%
%

%






\bibitem[\protect\citeauthoryear{Gnedenko}{Gnedenko \& Kolmogorov}{1968}]{GK1968}
Gnedenko, B.V. \& Kolmogorov, A.N. (1968)
\newblock {\em Limit Distributions for Sums of Independent Random Variables.}
\newblock Addison-Wesley.

\bibitem[\protect\citeauthoryear{Gregoire}{Gregoire}{1984}]{Gregoire1984}
Gregoire, G. (1984)
\newblock Negative binomial distributions for point processes.
\newblock {\em Stoch.  Proc. Appl.},  16, 179--188.

\bibitem[\protect\citeauthoryear{Griffiths}{Griffiths}{1979}]{Griffiths1979}
Griffiths, R.C. (1979)
\newblock On the distribution of allele frequencies in a diffusion model.
\newblock {\em Theoret. Pop. Biol.}, 15, 140--158.

\bibitem[\protect\citeauthoryear{Grote \& Speed}{Grote \& Speed}{2002}]{GSpeed2002}
Grote, M.N \& Speed, T.P. (2002).
\newblock Approximate Ewens formulae for symmetric overdominance selection.
\newblock {\it Ann. Appl. Probab.}, 12, 637--663.


\bibitem[\protect\citeauthoryear{Handa}{Handa}{2009}]{Handa2009}
Handa, K. (2009) 
\newblock The two-parameter Poisson-Dirichlet point process.
\newblock {\em Bernoulli}, 15, 1082--1116.


\bibitem[\protect\citeauthoryear{Hansen}{Hansen}{1994}]{Hansen1994}
Hansen, J.  (1994) 
\newblock Order statistics for decomposable combinatorial structures
\newblock {\em Random Structures \& Algorithms}, 5I, 517--533.

\bibitem[\protect\citeauthoryear{Hensley}{Hensley}{1982}]{Hensley1982}
Hensley, D. (1982) 
\newblock The convolution powers of the Dickman function.
\newblock {\em J. Lond. Math. Soc.}, s2-33, 395--406.


%
%
%
%


\bibitem[\protect\citeauthoryear{Ipsen \& Maller}{Ipsen \& Maller}{2017}]{IpsenMaller2017}
Ipsen, Y.F. \& Maller, R.A. (2017)
\newblock Negative binomial construction of random discrete distributions on   the infinite simplex.
\newblock {\em Theor. Stoch. Proc.}, 22, 34--46.



\bibitem[\protect\citeauthoryear{Ipsen, Maller \& Shemehsavar}{Ipsen, Maller \& Shemehsavar}{2020}]{IpsenMallerShemehsavar2020}
Ipsen, Y.F.,  Maller, R.A. \& Shemehsavar (2020)
\newblock  Limiting distributions of generalised Poisson-Dirichlet  distributions based on negative binomial processes.
\newblock {\em J. Theor. Prob.}, 33, 1974--2000.

\bibitem[\protect\citeauthoryear{Ipsen, Shemehsavar \& Maller}{Ipsen,  Shemehsavar \& Maller}{2018}]{IpsenMallerShemehsavar2018}
Ipsen, Y.F., Shemehsavar, S. \&  Maller, R.A.  (2018)
\newblock Species sampling models generated by negative binomial processes. 
\newblock  {\em arXiv: 1904.13046}.

\bibitem[\protect\citeauthoryear{Ipsen, Maller \& Shemehsavar}{Ipsen, Maller \& Shemehsavar}{2020}]{IpsenMallerShemehsavar2020}
Ipsen, Y.F.,  Maller, R.A. \& Shemehsavar, S. (2020)
\newblock  Size biased sampling from the Dickman subordinator.
\newblock {\em Stoch. Proc. Appl.}, 130, 6880--6900.

\bibitem[\protect\citeauthoryear{Ipsen, Maller \& Shemehsavar}{Ipsen, Maller \& Shemehsavar}{2021}]{IpsenMallerShemehsavar2021}
Ipsen, Y.F.,  Maller, R.A. \& Shemehsavar, S. (2021)
\newblock  A generalised Dickman distribution and the number of species in a negative binomial process model.
\newblock {\em Adv. Appl. Prob., 53, 370--399}.


\bibitem[\protect\citeauthoryear{LancelotJames2008}
{James}{2008}]{LancelotJames2008}
 James, Lancelot F.  (2008)
\newblock Large sample asymptotics for the
two-parameter Poisson–Dirichlet process.
\newblock {\em Pushing the Limits of Contemporary Statistics: Contributions in Honor of Jayanta K. Ghosh, Vol. 3 (2008) 187–199,  Inst. Math. Statistics, 2008}.

\bibitem[\protect\citeauthoryear{Joyce, Krone, \& Kurtz}
{Joyce, Krone, \& Kurtz}{2002}]
{JoyceKroneKurtz2002}
Joyce, P.,  Krone, S.M. \& Kurtz, T.G. (2002) 
\newblock  Gaussian limits associated with the Poisson-Dirichlet distribution and the Ewens sampling formula.
\newblock {\em Ann. Appl. Prob.}, 12, 101--124.

\bibitem[\protect\citeauthoryear{Kingman}{Kingman}{1975}]{kingman1975}
Kingman, J.F.C. (1975)
\newblock Random discrete distributions.
\newblock {\em J. Roy. Statist. Soc. B},  37, 1--22.


\bibitem[\protect\citeauthoryear{Lijoi, Mena \& Prunster}
{Lijoi, Mena \& Prunster}{2005}]
{LijoiMenaPrunster2005}
 Lijoi, A., Mena, R.H. \& Prunster, I.  (2005) 
\newblock  Mixture modeling with normalized inverse-Gaussian priors.
\newblock {\em J. Amer. Statist. Assoc.}, 100, 1278-1291.


\bibitem[\protect\citeauthoryear{Maller \& Shemehsavar}{Maller \& Shemehsavar}{2023}]{MallerShemehsavar2023}
Maller, R.A. \& Shemehsavar, S. (2023)
\newblock  Generalized Poisson--Dirichlet distributions based on the Dickman subordinator.
\newblock {\em Theor. Prob. Appl}, 67, 593-612.

%
%

%

%
%
%

%


%






%


%
\bibitem[\protect\citeauthoryear{Perman}{Perman}{1993}]{perman1993}
Perman, M. (1993).
\newblock Order statistics for jumps of normalised subordinators.
\newblock {\em Stoch.  Proc. Appl.},  46, 267--281.

\bibitem[\protect\citeauthoryear{Perman, Pitman \& Yor}{Perman, Pitman \& Yor}{1992}]{PPY1992}
Perman, M., Pitman, J., \& Yor, M. (1992)
\newblock Size-biased sampling of {P}oisson point processes and excursions.
\newblock {\em Probab. Theor.  Rel. Fields},  92, 21--39.

%

\bibitem[\protect\citeauthoryear{Pitman}{Pitman}{1995}]{Pitman1995PTRF}
Pitman, J. (1995)
\newblock Exchangeable and partially exchangeable random partitions.
\newblock {\em Probab. Theory Relat. Fields}, 102, 145--158.

%

\bibitem[\protect\citeauthoryear{Pitman}{Pitman}{1997}]{Pitman1997}
Pitman, J. (1997)
\newblock Partition structures derived from Brownian
 motion and stable subordinators.
\newblock {\em Bernoulli}, 3, 79--96.


\bibitem[\protect\citeauthoryear{Pitman}{Pitman}{2006}]{Pitman2006}
Pitman, J. (2006)
\newblock {\em Combinatorial Stochastic Processes}.
\newblock Springer-Verlag, Berlin.


\bibitem[\protect\citeauthoryear{Pitman \& Yor}{Pitman \& Yor}{1997}]{PY1997}
Pitman, J. \& Yor, M. (1997)
\newblock The two-parameter {Poisson--Dirichlet} distribution derived from a
  stable subordinator.
\newblock {\em Ann. Probab.}, 25, 855--900.
%



\bibitem[\protect\citeauthoryear{Ruggiero, Walker \& Favaro}{Ruggiero, Walker \& Favaro}{2013}]{RWF2013}
Ruggiero, M., Walker, S.G. \& Favaro, S. (2013)
\newblock Alpha-diversity processes and normalized
 inverse-Gaussian diffusions.
\newblock {\em Ann. Appl. Probab.}, 23,  386--425.


\bibitem[\protect\citeauthoryear{Saa \& Venegeroles}{Saa \& Venegeroles}{2011}]{SaaVen2011}
 Saa, A. \& Venegeroles, R. (2011)
\newblock  Alternative numerical computation of one-sided L\'evy and Mittag-Leffler distributions.
\newblock {\em Physical Review E}, 84 (2 Pt 2):026702.




%




\bibitem[\protect\citeauthoryear{Watterson}{Watterson}{1977}]{Watterson1974}
Watterson, G.A. (1974)
\newblock Models for the logarithmic species abundance distribution
\newblock {\em Theor. Pop. Biol.}, 6, 217--250.


%


\bibitem[\protect\citeauthoryear{Zhou, Favaro \& Walker}{Zhou, Favaro \& Walker}{2017}]{ZFW2017}
Zhou, M., Favaro, S. \& Walker, S.G. (2017)
\newblock 
Frequency of frequencies distributions and size-dependent exchangeable 
random partitions
\newblock {\it J. Amer. Statist. Assoc.}, 112, 1623--1635.

\end{thebibliography}
\end{document}